\documentclass[reqno]{amsart}
\usepackage{amssymb}
\usepackage{amsmath}

\makeatletter
\@addtoreset{equation}{section}
\makeatother

\renewcommand\thefigure{\thesection.\@arabic\c@figure}
\renewcommand\thetable{\thesection.\@arabic\c@table}

\newtheorem{theorem}{Theorem}[section]
\newtheorem{lemma}[theorem]{Lemma}
\newtheorem{proposition}[theorem]{Proposition}
\newtheorem{corollary}[theorem]{Corollary}

\newcommand{\mc}[1]{{\mathcal #1}}
\newcommand{\mf}[1]{{\mathfrak #1}}
\newcommand{\mb}[1]{{\mathbf #1}}
\newcommand{\bb}[1]{{\mathbb #1}}

\newcommand{\V}{|\! |\! |}

\def\emptysquare{{\hbox{\vrule height6pt width0.6pt depth0pt%
\vbox{\hrule height0.6pt width4.8pt depth0pt%
\vglue4.8pt%
\hrule height0.6pt width4.8pt depth0pt}%
\vrule height6pt width0.6pt depth0pt}}}

\def\qed{\unskip\nobreak
\hfil\penalty50\hskip1.75em\null\nobreak\hfil\emptysquare
{\parfillskip=0pt \finalhyphendemerits=0 \par}\medskip}

\begin{document}

\title[asymmetric exclusion processes under diffusive
scaling]{Hydrodynamic limit of asymmetric exclusion processes under
  diffusive scaling in $d\ge 3$.}

\author{C. Landim, M. Sued and G. Valle}

\address{\noindent IMPA, Estrada Dona Castorina 110,
CEP 22460 Rio de Janeiro, Brasil and CNRS UMR 6085,
Universit\'e de Rouen, 76128 Mont Saint Aignan, France.
\newline
e-mail:  \rm \texttt{landim@impa.br}
}

\address{\noindent IMPA, Estrada Dona Castorina 110,
CEP 22460 Rio de Janeiro, Brasil.
\newline
e-mail:  \rm \texttt{msued@impa.br}
}

\address{\noindent IMPA, Estrada Dona Castorina 110,
CEP 22460 Rio de Janeiro, Brasil.
\newline
e-mail:  \rm \texttt{valle@impa.br}
}

\dedicatory{Dedicated to J\'ozsef Fritz on his sixtieth birthday.}

\begin{abstract}
  We consider the asymmetric exclusion process. We start from a
  profile which is constant along the drift direction and prove that
  the density profile, under a diffusive rescaling of time, converges
  to the solution of a parabolic equation.
\end{abstract}

\maketitle

\section{Introduction}
\label{sec1}

Consider the asymmetric exclusion process evolving on the lattice $\bb
Z^d$.  This dynamics can be informally described as follows~: fix a
translation invariant transition probability $p(x,y) = p(0,y-x)
=p(y-x)$. Each particle, independently from the others, waits a mean
one exponential time, at the end of which being at $x$ it chooses the
site $x+y$ with probability $p(y)$. If the chosen site is unoccupied,
the particle jumps, otherwise it stays where it is. In both cases,
after its attempt, the particle waits a new mean one exponential time.

The configurations of the state space $\{0,1\}^{\bb Z^d}$ are denoted
by the Greek letter $\eta$ so that, for $x$ in $\bb Z^d$, $\eta (x)$
is equal to $1$ or $0$, whether site $x$ is occupied or not.  For each
density $0\le \alpha\le 1$, the Bernoulli product measure with
parameter $\alpha$, denoted by $\nu_\alpha$, is invariant.

The macroscopic evolution of the process under Euler rescaling is
described \cite{r} by the first order quasilinear hyperbolic equation
\begin{equation}
\label{f1}
\partial_t \rho \; +\; q \cdot \nabla  F(\rho) \; =\; 0\; ,
\end{equation}
where $F(a) = a(1-a)$ and $q \in\bb R^d$ is the mean drift of each
particle~: $q =\sum_z z p(z)$. Assume that the system starts from a
product measure with slowly varying density $\rho_0( \varepsilon u)$.
Under Euler scaling (times of order $t\varepsilon^{-1}$) the density
has still a slowly varying profile $\lambda_\varepsilon (t,
\varepsilon u)$ which converges weakly (in fact pointwisely at every
continuity point, \cite{l}) to the entropy solution of equation
(\ref{f1}) with initial data $\rho_0$.

In the context of asymmetric interacting particle systems the
Navier--Stokes equations takes the form
\begin{equation}
\label{f2}
\partial_t \rho^\varepsilon \; 
+\; q \cdot \nabla F(\rho^\varepsilon) \; =\; \varepsilon
\sum_{i,j} \partial_{u_i} \Big( a_{i,j}(\rho^\varepsilon) \partial_{u_j}
\rho^\varepsilon \Big )\; ,
\end{equation}
where $a$ is a diffusion coefficient.  Three different interpretations
have been proposed for the Navier--Stokes corrections~:

\noindent  {\bf (a) The incompressible limit} (\cite{em},
\cite{emy1})~: Consider a small perturbation of a constant profile
$\alpha_0$~: $\rho^\varepsilon_0 = \alpha_0 + \varepsilon \varphi$.
Assuming that this form persists at latter times ($\rho^\varepsilon
(t,u) = \alpha_0 + \varepsilon \varphi (t,u)$) we obtain from
(\ref{f2}) the following equation for $\varphi_\varepsilon =
\varphi (t\varepsilon^{-1},u)$
$$
\partial_t \varphi _\varepsilon \; +\; \varepsilon^{-1}
F'(\alpha_0) q \cdot \varphi_\varepsilon \; +\; (1/2)
F''(\alpha_0)q \cdot \nabla \varphi_\varepsilon^2 \; =\;
a_{i,j}(\alpha_0) \sum_{i,j} \partial_{u_i, u_j}^2 \varphi_\varepsilon
\; +\; O(\varepsilon)\; .
$$
A Galilean transformation $m_\varepsilon(t,u) =
\varphi_\varepsilon(t,u + \varepsilon^{-1} t F'(\alpha_0) q)$
permits to remove the diverging term of the last differential equation
and to get a limit equation for $m = \lim_{\varepsilon\to 0}
m_\varepsilon$
$$
\partial_t m \; +\; (1/2) F''(\alpha_0) q \cdot \nabla m^2 \;
=\; a_{i,j}(\alpha_0) \sum_{i,j} \partial_{u_i, u_j}^2 m \; .
$$

\noindent {\bf (b) First order correction to the hydrodynamic equation} 
(\cite{d}, \cite{loy1})~: Fix a smooth profile $\rho_0\colon\bb
R^d\to\bb R_+$ and consider a process starting from a product measure
with slowly varying density $\rho_0(\varepsilon u)$. We have seen that
under Euler scaling the density is still a slowly varying profile
$\lambda_\varepsilon (t, \varepsilon u)$ which converges weakly to the
entropy solution of equation (\ref{f1}) with initial data $\rho_0$.
This second interpretation asserts that the solution of equation
(\ref{f2}) with initial profile $\rho_0$ approximates
$\lambda_\varepsilon$ up to the order $\varepsilon$~:
$$
\varepsilon^{-1} \big(\lambda_\varepsilon -  \rho^\varepsilon \big)\to \; 0
$$
in a weak sense as $\varepsilon \downarrow 0$.

\noindent{\bf (c) Long time behaviour} (\cite{loy1}, \cite{bkl})~: The third 
interpretation consists in analyzing the behaviour of the solution of
equation (\ref{f2}) in time scales of order $t\varepsilon^{-1}$. Let
$b_\varepsilon (t,u) = \rho(t\varepsilon^{-1}, u)$. From (\ref{f2}) we
obtain the following equation for $b_\varepsilon$~:
$$
\partial_t b_\varepsilon \; +\; \varepsilon^{-1} q \cdot 
\nabla F(b_\varepsilon) \; =\; 
\sum_{i,j} \partial_{u_i} \Big( a_{i,j}(b_\varepsilon) \partial_{u_j}
b_\varepsilon \Big)\; .
$$
To eliminate the diverging term $\varepsilon^{-1} q \cdot \nabla
F(b_\varepsilon)$, assume that the initial data (and therefore the
solution at any fixed time) is constant along the drift direction~: $q
\cdot \nabla\rho_0 =0$. In this case we get the parabolic equation
$$
\partial_t b_\varepsilon  \; =\; 
\sum_{i,j} \partial_{u_i} \Big( a_{i,j}(b_\varepsilon) \partial_{u_j}
b_\varepsilon \Big)
$$
which describes the evolution of the system in the hyperplane
orthogonal to the drift.
\smallskip

Notice that while the first and the third interpretation concern the
behaviour of the system under diffusive rescaling, the second one is a
statement on the process under Euler rescaling.  Interpretation (a)
and (b) have been proved \cite{emy1}, \cite{loy1} for asymmetric
simple exclusion processes in dimensions $d\ge 3$ and a double
variational formula for the diffusion coefficient was deduced.  As one
would expect, the diffusion coefficients of the two interpretations
are the same and may be expressed by a Green--Kubo formula \cite{ly}.
It was also proved (Corollary 6.2, \cite{loy2}) that the diffusion
coefficient is strictly bounded below in the matrix sense by the
diffusion coefficient that governs the evolution of the symmetric
process and that it depends smoothly on the density \cite{lov2}.

In contrast with interpretations (a) and (b), the third one is
meaningful in dimension $d\ge 2$. It has been proved in \cite{bkl} for
asymmetric zero range processes. The purpose of this paper is to give
a rigorous proof of the third interpretation in the case of asymmetric
exclusion processes in dimension $d\ge 3$. The proof in this context
is much more demanding because the process is nongradient.  In
particular, we obtain a non-trivial diffusion coefficient.

\section{Notation and Results}
\label{sec2}

Fix a finite range probability measure $p(\cdot)$ on $\bb Z^d$.  The
exclusion process evolving on the discrete torus $\bb T_N^d = \{0,
\dots, N-1\}^d$ associated to $p(\cdot)$ is the Markov process on the
state space $\mc X_N=\{0,1\}^{\bb T_N^d}$ whose generator $L_N$ acts
on a local function $f$ as
\begin{equation}
\label{eq:1}
(L_N f)(\eta)\; =\;  \sum_{\scriptstyle x,y \in\bb T_N^d}
p(y)  \eta(x) \{ 1-\eta(x+y) \} 
[f( \sigma^{x,x+y} \eta)-f(\eta)] \;,
\end{equation}
where $\sigma^{x,x+y} \eta$ is the configuration obtained from $\eta$ by
exchanging the occupation variables $\eta(x)$, $\eta(x+y)$:
$$
(\sigma^{x,x+y} \eta)(z) \;=\;
\left\{
\begin{array}{ll}
\eta(z) & \hbox{if $z\not = x$, $x+y$\; ,} \\
\eta(x) & \hbox{if $z= x+y$\; ,} \\
\eta(x+y) & \hbox{if $z= x$\; .} \\
\end{array}
\right.
$$

Fix $\alpha$ in $(0,1)$ and denote by $\nu_\alpha^N$ the Bernoulli
product measure on $\mc X_N$ with density $\alpha$.  Let $L_N^*$ be
the adjoint of $L_N$ in $L^2(\nu_\alpha^N)$.  This operator is
obtained by replacing $p(y)$ by $p^*(y) = p(-y)$ in (\ref{eq:1}).

Denote by $\bb T^d$ the $d$-dimensional torus. Fix a continuous
function $\rho_0: \bb T^d\to [0,1]$ and denote by
$\nu^N_{\rho_0(\cdot)}$ the product measure on $\{0,1\}^{\bb T_N^d}$
associated to $\rho_0$. This is the Bernoulli product measure on
$\{0,1\}^{\bb T_N^d}$ with marginals given by
$$
\nu^N_{\rho_0(\cdot)} \{ \eta (x) =1 \} \; =\; \rho_0(x/N)
$$
for $x$ in $\bb T_N^d$. 

For $N\ge 1$ and a configuration $\eta$, denote by $\pi^N (\eta)$ the
empirical measure associated to $\eta$. This is the measure on $\bb
T^d$ obtained by assigning mass $N^{-d}$ to each particle of $\eta$:
$$
\pi^N (\eta) \;=\; N^{-d} \sum_{x\in \bb T_N^d} \eta(x) \delta_{x/N}
\;, 
$$
where $\delta_u$ stands for the Dirac measure on $u$. It has been
proved in \cite{r} that if particles are initially distributed
according to $\nu^N_{\rho_0(\cdot)}$ for some profile $\rho_0: \bb
T^d\to [0,1]$, then $\pi^N(\eta_{tN})$ converges in probability to
$\rho (t,u) du$, where $\rho$ is the entropy solution of the Burgers
equation 
\begin{equation}
\label{eq:2}
\partial_t \rho \; +\; q \cdot \nabla  F(\rho) \; =\; 0\; ,
\end{equation}
where $F(a) = a(1-a)$ and $q \in\bb R^d$ is the mean drift of each
particle: $q =\sum_z z p(z)$.

In this article, we investigate the diffusive behavior of the
empirical measure $\pi^N$, that is, its evolution in times of order
$N^2$. 

As time increases, the solution of Burgers equation (\ref{f2})
converges to a stationary profile which is constant along the drift
direction: 
$$
\lim_{t\to\infty} \rho(t,u) \; =\; \rho_\infty(u) \; =\; \int_0^1
\rho_0 (u + r q)\, dr\; ,
$$
provided $\rho_0$ stands for the initial data. The limit should be
understood pointwisely. In particular, in a time scale of order $N^2$,
the profile of the empirical measure should immediately become constant
along the drift direction.

We shall therefore assume that the initial state is a product measure
$\nu^N_{\rho_0(\cdot)}$ associated to a profile $\rho_0$ constant
along the drift direction:
\begin{equation}
\label{eq:3}
q \cdot \nabla \rho_0 (u) \;=\; 0
\end{equation}
for all $u$ in $\bb T^d$. Assume furthermore that the profile is
bounded away from $0$ and $1$:
\begin{equation}
\label{eq:5}
\delta_0 \;\le\; \rho_0(u) \;\le\; 1- \delta_0
\end{equation}
for some $\delta_0>0$.  

\begin{theorem}
\label{s1}
Assume that the initial state is distributed according to
$\nu^N_{\rho_0(\cdot)}$, where the profile $\rho_0$ satisfies
(\ref{eq:3}), (\ref{eq:5}). There exists a smooth matrix-valued
function $a(\alpha) = \{ a_{i,j} (\alpha), \, 1\le i,j\le d\}$ with
the following property. For each $t\ge 0$, $\pi^N(\eta_{tN^2})$
converges in probability to $\rho (t,u) du$, where the density $\rho$
is the solution of the parabolic equation
\begin{equation}
\label{eq:4}
\left\{
\begin{array}{l}
\partial_t\rho \;=\; \sum_{i,j} \partial_{u_i}(a_{i,j}(\rho)
\partial_{u_j}\rho)\;,\\
\rho(0,\cdot) \;=\; \rho_0(\cdot)\;.
\end{array}
\right.
\end{equation}
\end{theorem}

In this theorem, $a_{i,j}(\alpha) = D_{i,j}(\alpha) +
(1/2)(1-2\alpha) \sigma_{i,j}$, where $D_{i,j}(\alpha)$ is the matrix
given by (\ref{eq:27}) and $\sigma_{i,j}$ the covariance matrix of the
transition probability $p(\cdot)$:
$$
\sigma_{i,j} \;=\; \sum_{y\in\bb Z^d} p(y) \, y_i\, y_j\;.
$$
Notice that by the maximum principle, $\delta_0\le \rho(t,u) \le
1-\delta_0$ for all $(t,u)$. Moreover, the solution of the
hydrodynamic equation is constant along the drift direction,
$$
\sum_{i=1}^d q_i \, (\partial_{u_i} \rho)(t,u) \;=\;0
$$
because so is the initial data.

This theorem is an elementary consequence of the following estimate on
the relative entropy of the state of the process with respect to a
local Gibbs state. For two measures $\mu$, $\nu$ on $\{0,1\}^{\bb
  T_N^d}$, denote by $H_N(\mu | \nu)$ the relative entropy of $\mu$ with
respect to $\nu$:
$$
H_N(\mu \,|\, \nu) \;=\; \sup_{f} \Big\{ \int f \, d\mu \;-\; \log \int
e^f d\nu \Big\}\;,
$$
where the supremum is carried over all bounded, continuous
functions, which in our finite setting coincide with all functions.
For $t\ge 0$, denote by $S^N_t$ the semigroup associated to the Markov
process with generator (\ref{eq:1}) speeded up by $N^2$.

\begin{theorem}
\label{s2}
Under the assumptions of Theorem \ref{s1} on the initial profile
$\rho_0$, let $\{ \mu_N, \, N\ge 1 \}$ be a sequence of probability
measures on $\{0,1\}^{\bb T_N^d}$ whose entropy with respect to
$\nu^N_{\rho_0(\cdot)}$ is of order $o(N^d)$:
$$
H_N(\mu_N \,|\,\nu^N_{\rho_0(\cdot)}) \;=\; o(N^d)\;.
$$
Then, for every $t\geq 0$, the relative entropy of the state of the
process at time $tN^2$ with respect to $\nu^N_{\rho(t,\cdot)}$ is also
of order $o(N^d)$:
$$
H_N (\mu_N S^N_t \,|\, \nu^N_{\rho(t,\cdot)}) \;=\; o(N^d)\;,
$$
provided $\rho(t,u)$ is the solution of (\ref{eq:4}).
\end{theorem}

In view of this result, we can weaken the assumptions of Theorem
\ref{s1} and assume only that the initial state has relative entropy
of order $o(N^d)$ with respect to $\nu^N_{\rho_0(\cdot)}$.

\section{Relative entropy estimates}
\label{sec3}

We introduce in this section some auxiliary measures which will play a
central role in the proof of Theorem \ref{s2}. The statements
presented here appeared essentially in the same form in \cite{emy1}
and \cite{loy1}. We include their proof in sake of completeness.

Fix a profile $\rho_0$ constant along the drift direction and bounded
away from $0$ and $1$ as in (\ref{eq:5}). Denote by $\rho(t,u)$ the
smooth solution of the parabolic equation (\ref{eq:4}). Fix
$0<\alpha<1$. For $N\ge 1$, denote by $f^N_t$ the density of $\mu^N
S^N_t$ with respect to $\nu_\alpha^N$. An elementary computation shows
that $f^N_t$ is the solution of
$$
\partial_t f^N_t \;=\; N^2 L_N^* f^N_t\; ,
$$
where $L_N^*$ is the adjoint of the generator $L_N$ in
$L^2(\nu_\alpha^N)$. 

Denote by $\mf F$ the space of functions $\mf f \colon [0,1]\times
\{0,1\}^{\bb Z ^d} \to\bb R$ such that

\begin{enumerate}
\item There exists a finite set $\Lambda$ such that for each $\beta$
  in $[0,1]$ the support of $\mf f(\beta,\cdot)$ is contained in
  $\Lambda$.

\smallskip
\item For each configuration $\eta$, $\mf f (\cdot,\eta)$ is a smooth
  function.

\smallskip
\item For each density $\beta$, the cylinder functions $\mf f(\beta,
  \cdot)$, $\mf f_1 (\beta, \cdot)$ have zero mean with respect to
  $\nu_\beta$. Here, $\mf f_1 (\beta, \cdot)$ stands for the
  derivative of $\mf f(\beta, \eta)$ with respect to the first
  coordinate.
\end{enumerate}

Let $\lambda : \bb R_+\times \bb T^d\to \bb R$ be defined by 
$$
\lambda (t, u)\; =\; \log \frac{\rho (t,u) (1-\alpha) }
{\alpha [ 1- \rho (t,u) ]} \;\cdot
$$
$\lambda (t,u)$ is well defined because the solution $\rho(t,u)$ of
the hydrodynamic equation (\ref{eq:4}) is bounded away from $0$ and
$1$.  Denote by $\psi^{N}_{t} (\eta)$ the density of $\nu^N_{\rho(t,
  \cdot)}$ with respect to $\nu_\alpha$:
$$
\psi^{N}_{t} (\eta) \;=\; \frac{1}{Z_t} \, 
\exp\Big \{ \sum_{x\in\bb T_N^d} \lambda (t, x/N) \eta (x) \Big\} \; ,
$$
where $Z_t$ is a renormalizing constant.

For functions $\mf f_i$ in $\mf F$, $1\le i\le d$, a time $t\ge 0$ and
integers $\ell \ll M \ll N$, define the density $\psi^N_{t,\mf f}
(\eta)= \psi^{N,M,\ell}_{t,\mf f} (\eta)$ with respect to the reference
measure $\nu_\alpha^N$ by
\begin{eqnarray*}
\!\!\!\!\!\!\!\!\!\!\! &&
\psi^N_{t,\mf f}(\eta)\;  =\; \frac{1}{Z_t^{\mf f}} \, 
\exp\Big \{ \sum_{x\in\bb T_N^d} \lambda (t, x/N) \eta (x) \Big\} 
\;\times \\
\!\!\!\!\!\!\!\!\!\!\! && \qquad\qquad\qquad \times\;
\exp\Big \{ -\; N^{-1} \sum_{i=1}^d \sum_{x\in\bb T_N^d} 
\partial_{u_i} \lambda  (t,x/N) \frac 1{|\Lambda_{\ell'}|} 
\sum_{y\in \Lambda_{\ell'}} \mf f_i (\eta^M (x) ,\tau_{x+y} \eta ) 
\Big\}\; ,
\end{eqnarray*}
where $Z_t^{\mf f}$ is a renormalizing constant, $\Lambda_K = \{-K,
\dots , K\}^d$ is a cube of length $2K+1$ centered at the origin,
$\eta^K (x)$ is the mean density of particles in $x+\Lambda_K$:
\begin{equation}
\label{eq:14}
\eta^K (x) \;=\; \frac 1{|\Lambda_K|} \sum_{y\in x + \Lambda_K}
\eta(y) 
\end{equation}
and $\ell'=\ell -A$ for a finite constant $A$ chosen for the support
of $\mf f_i (\beta, \tau_y\eta)$ to be contained in $\Lambda_\ell$ for
all $1\le i\le d$, $|y|\le \ell - A$. Throughout this article, $A$
stands for a finite integer related to the support of the transition
probability $p(\cdot)$ or to the support of some local function.

In the following, we will need to take $M$ as a function of $N$ and
$\ell$ as an independent integer which increases to $\infty$ after
$N$. In fact we will require $M$ to be such that
\begin{equation}
\label{eq:16}
\lim_{N\to\infty} \frac{|\Lambda_M|} N \;=\; 0\;, \quad
\lim_{N\to\infty} \frac N { M |\Lambda_M|} \;=\; 0\;. 
\end{equation}

We present three elementary results which illustrate some properties
of the density $\psi^N_{t,\mf f}(\eta)$. Denote by $s_{\mf f}$ the
smallest integer $m$ with the property that the common support of
the local functions $\mf f_i (\beta , \cdot)$, $1\le i\le d$,
$0\le\beta\le 1$, is contained in $\Lambda_m$.

\begin{lemma}
\label{s9}
Assume that $s_{\mf f} \le \ell \le M$ and that $\lim_{N\to\infty}
|\Lambda_M|/ N = 0$. Fix a density $f$ with respect to the reference
measure $\nu_\alpha^N$. There exists a finite constant $C$, depending
only on $\mf f$ and $\rho(t,u)$, such that
$$
\big\vert\, 
H_N (f \,|\, \nu^N_{\rho(t, \cdot)}) \;-\;
H_N (f \,|\, \psi^N_{t,\mf f} ) \, \big\vert\,  \;\le\; C N^{d-1}
$$
for all $N\ge 1$.
\end{lemma}

In the statement of this result and frequently in this article, if
measures $\mu$, $\nu$ have density $f$, $g$ with respect to the
reference measure $\nu_\alpha^N$, to keep notation simple, we denote
by $H_N (f \,|\, g )$ the entropy of $f \, d\nu_\alpha^N$ with respect
to $g \, d\nu_\alpha^N$ and by $E_{f} [\cdot]$ the expectation with
respect to $f \, d\nu_\alpha^N$.

\begin{proof}
Fix a density $f$.  By the explicit formula for the entropy, the
difference $H_N (f \,|\, \nu^N_{\rho(t, \cdot)}) - H_N (f \,|\,
\psi^N_{t,\mf f} )$ is equal to
$$
\int f \, \log \frac{\psi^N_{t,\mf f}}{\psi^N_t}\, d\nu_\alpha^N
\;=\; O(N^{d-1}) \;-\; \log \frac{Z_t^{\mf f}}{Z_t}\; \cdot
$$
In particular, we just need to show that the second term on the
right hand side is absolutely bounded by $C N^{d-1}$. By definition of
the renormalizing constant $Z_t^{\mf f}$, $Z_t$, the logarithm is
equal to
\begin{equation}
\label{eq:15}
\log E_{\nu^N_{\rho(t, \cdot)}} \Big[ \exp\Big\{ 
-\; N^{-1} \sum_{i=1}^d \sum_{x\in\bb T_N^d} 
\partial_{u_i} \lambda  (t,x/N) (A_\ell \mf f_i) 
(\eta^M (x) ,\tau_{x} \eta ) \Big\} \Big] \;,
\end{equation}
where, for a function $\mf f$ in $\mf F$ and a positive integer $\ell$,
$$
(A_\ell \mf f)( \beta, \eta) \;=\;
\frac 1{|\Lambda_{\ell'}|} \sum_{y\in \Lambda_{\ell'}} 
\mf f (\beta ,\tau_{y} \eta )\; .
$$
By Jensen inequality, (\ref{eq:15}) is bounded below by $-C
N^{d-1}$. On the other hand, since $\ell \le M$, $A_\ell \mf f_i
(\eta^M (0) ,\eta )$ depends on the configuration $\eta$ only through
$\eta (z)$ for $z$ in $\Lambda_M$. In particular, since
$\nu^N_{\rho(t, \cdot)}$ is a product measure, by H\"older inequality,
(\ref{eq:15}) is bounded above by
$$
\frac 1{ d |\Lambda_M|} \sum_{i=1}^d \sum_{x\in\bb T_N^d} 
\log E_{\nu^N_{\rho(t, \cdot)}} \Big[ \exp\Big\{ - d |\Lambda_M| N^{-1}
(\partial_{u_i} \lambda)  (t,x/N) (A_\ell \mf f_i) 
(\eta^M (x) ,\tau_{x} \eta ) \Big\} \Big]\; .
$$
Since by assumption $|\Lambda_M| N^{-1}$ vanishes as
$N\uparrow\infty$, we may expand the exponential up to the second
order to show that this expression is less than or equal to $C
N^{d-1}$. This concludes the proof of the lemma.
\end{proof}

Taking $f = \psi^N_{t,\mf f}$ in Lemma \ref{s9}, we obtain a bound on the
entropy of $\psi^N_{t,\mf f}$ with respect to $\nu^N_{\rho(t, \cdot)}$.

\begin{corollary}
\label{s10}
Under the assumptions of Lemma \ref{s9}, there exists a finite
constant $C$, depending only on $\mf f$ and $\rho(t,u)$, such that
$$
H_N (\psi^N_{t,\mf f} \,|\, \nu^N_{\rho(t, \cdot)})  
\;\le\; C N^{d-1}
$$
for all $N\ge 1$.
\end{corollary}

\begin{corollary}
\label{s11}
Fix a smooth function $H: \bb T^d \to R$ and a function $\mf g$ in
$\mf F$. There exists a finite constant $C_0$, depending only on $\mf
f$, $\mf g$, $H$ and $\rho(t,u)$, such that
\begin{eqnarray*}
\!\!\!\!\!\!\!\!\!\!\! &&
\Big\vert\,
E_{\psi^N_{t,\mf f}} \Big[ N^{-d} \sum_{x\in\bb T_N^d} 
H  (x/N) (A_\ell \mf g) (\eta^M (x) ,\tau_{x} \eta ) \Big] \;-\; \\
\!\!\!\!\!\!\!\!\!\!\! && \qquad\qquad\qquad
E_{\nu^N_{\rho(t, \cdot)}} \Big[ N^{-d} \sum_{x\in\bb T_N^d} 
H  (x/N) (A_\ell \mf g) (\eta^M (x) ,\tau_{x} \eta )\Big] 
\, \Big\vert
\;\le\; C_0 \sqrt{|\Lambda_M|/N}\;.
\end{eqnarray*}
\end{corollary}

\begin{proof}
By the entropy inequality
$$
E_{\psi^N_{t,\mf f}} \Big[ N^{-d} \sum_{x\in\bb T_N^d} 
H  (x/N) (A_\ell \mf g) (\eta^M (x) ,\tau_{x} \eta ) \Big]
$$
is less than or equal to
\begin{eqnarray*}
\!\!\!\!\!\!\!\!\!\!\! &&
\frac{H_N (\psi^N_{t,\mf f} \,|\, \nu^N_{\rho(t, \cdot)})}{\gamma N^{d-1}}
\\
\!\!\!\!\!\!\!\!\!\!\! && \quad
\;+\; \frac 1{\gamma N^{d-1}} \log 
E_{\nu^N_{\rho(t, \cdot)}} \Big[ \exp\Big\{ \gamma N^{-1}
\sum_{x\in\bb T_N^d} 
H  (x/N) (A_\ell \mf g) (\eta^M (x) ,\tau_{x} \eta ) \Big\}\Big]
\end{eqnarray*}
for every $\gamma>0$. By Corollary \ref{s10}, the first term is
bounded by $C \gamma^{-1}$. On the other hand, repeating the argument
presented in the proof of Lemma \ref{s9}, we show that the second term
is less than or equal to
\begin{eqnarray*}
\!\!\!\!\!\!\!\!\!\!\! &&
E_{\nu^N_{\rho(t, \cdot)}} \Big[ N^{-d} \sum_{x\in\bb T_N^d} 
H  (x/N) (A_\ell \mf g) (\eta^M (x) ,\tau_{x} \eta )\Big] \\
\!\!\!\!\!\!\!\!\!\!\! && \quad +\;
\frac {C \gamma |\Lambda_M|} {N} 
E_{\nu^N_{\rho(t, \cdot)}} \Big[ N^{-d} \sum_{x\in\bb T_N^d} 
H  (x/N)^2 (A_\ell \mf g) (\eta^M (x) ,\tau_{x} \eta )^2\Big]
\end{eqnarray*}
provided that $\gamma |\Lambda_M| N^{-1}$ vanishes as
$N\uparrow\infty$. In this formula, $C$ is a finite constant which
depends on $\mf g$ and $H$.  In particular, the difference appearing
inside the absolute value in the statement of the corollary is less
than or equal to
$$
\frac{C}{\gamma}
\;+\; \frac{C \gamma |\Lambda_M|} {N }\;\cdot
$$
Taking $\gamma = \sqrt{N/|\Lambda_M|}$, we show that this
expression is bounded by $C \sqrt{|\Lambda_M| / N}$. Replacing $H$ by
$-H$, we we conclude the proof of the corollary.
\end{proof}

\section{Proof of Theorem \ref{s2}}
\label{sec5}

We prove in this section Theorem \ref{s2}. In view of Lemma
\ref{s9}, Theorem \ref{s2} is a consequence of the following
result.

\begin{proposition}
\label{s3}
Fix a measure $\mu^N$ such that $H_N(\mu^N \,|\,
\nu_{\rho_0(\cdot)}^N) = o(N^d)$. Assume that the profile $\rho_0$
satisfies (\ref{eq:3}), (\ref{eq:5}). There exist sequences $\{\mf
f_{i,n} , n\ge 1\}$, $1\le i\le d$, of functions in $\mf F$ such that 
$$
\lim_{n\to\infty} \limsup_{\ell \to\infty} \limsup_{N\to\infty} N^{-d}
H_N ( \mu^N S_t^N \,\vert\, \psi^N_{t,\mf f_n} ) \;=\; 0
$$
for every $t\ge 0$. In this formula, $\mf f_n = (\mf f_{1,n}, \dots,
\mf f_{d,n})$. 
\end{proposition}

The proof of Proposition \ref{s3} is divided in several steps.  To
keep notation simple, denote by $H^{\mf f}_N (t)$ the relative entropy
of $\mu^N S_t^N$ with respect to $\psi^N_{t,\mf f}\, d\nu_\alpha^N$:
$$
H^{\mf f}_N(t) \;=\; H_N (\mu^N S_t^N \,\vert\, \psi^N_{t,\mf f} )
\;.
$$
In view of Lemma \ref{s9} and of Gronwall inequality, it is enough to
show that for every $t\ge 0$,
\begin{equation}
\label{eq:6}
H^{\mf f}_N(t) \;\le\; o(N^d,\mf f) \;+\; \gamma^{-1} \int_0^t ds\, 
H_N ( \mu^N S_s^N  \,\vert\, \nu_{\rho(s, \cdot)}^N )
\end{equation}
for some $\gamma>0$. Here, $o(N^d,\mf f)$ stands for a finite constant
such that 
$$
\lim_{n \to\infty} \limsup_{\ell \to\infty} \limsup_{N\to\infty}
N^{-d} o(N^d, \mf f_n) \;=\; 0\;.
$$

The sequence $\{\mf f_{i,n},\, n\ge 1\}$ is given by Theorem \ref{s7}.
To keep notation simple, we perform all computations for a single
function $\mf f =(\mf f_1, \dots, \mf f_d)$ and then replace it by the
sequence $\mf f_n$.

Recall that $M$ depends on $N$ through the relations (\ref{eq:16}) and
that $\ell$ is an integer independent of $N$ which increases to
infinity after $N$. To prove (\ref{eq:6}), we start computing the
time derivative of the entropy $H^{\mf f}_N(t)$.  On the one hand, a
celebrated estimate of \cite{y} gives that
\begin{equation}
\label{eq:7}
\frac d{dt} H^{\mf f}_N (t) \;\le\; \int f^N_t \Big\{
\frac{N^2 L_N^* \psi^N_{t,\mf f}}{\psi^N_{t,\mf f}}-\partial_t
\log(\psi^N_{t,\mf f})\Big\}\, d\nu_\alpha^N\;.
\end{equation}
On the other hand, a straightforward computation, presented in section
\ref{sec4}, shows that the expression inside braces in the previous
integral is equal to
\begin{eqnarray}
\label{eq:8}
\!\!\!\!\!\!\!\!\!\!\!\! && 
N \sum_{i=1}^d \sum_{x\in\bb T_N^d} (\partial_{u_i} \lambda) (t, x/N)
\Big\{ \tau_x W_i^* - L_N^* (A_\ell \mf f_i) (\eta^M(x), \tau_x \eta)
\Big\}  \\
\!\!\!\!\!\!\!\!\!\!\!\! && 
\quad +\; (1/2) \sum_{i, j =1}^d \sum_{x\in\bb T_N^d} 
(\partial^2_{u_i, u_j} \lambda) (t, x/N) \, \tau_x G_{i,j} (\eta)
\nonumber \\
\!\!\!\!\!\!\!\!\!\!\!\! && 
\qquad +\; (1/2) \sum_{i,j=1}^d \sum_{x\in\bb T_N^d} (\partial_{u_i}
\lambda) (t, x/N) (\partial_{u_j} \lambda) (t, x/N) \, \tau_x
H_{i,j}(\eta) 
\nonumber \\
\!\!\!\!\!\!\!\!\!\!\!\! && 
\quad \qquad -\; \sum_{x\in\bb T_N^d} (\partial_{t} \lambda) 
(t, x/N) \eta (x) \; +\; E_{\psi^N_{t,\mf f}} \Big[ 
\sum_{x\in\bb T_N^d} (\partial_{t} \lambda) (t, x/N) \eta (x) \Big]
\;+\; o(N^d)\;. \nonumber
\end{eqnarray}
In this formula, $o(N^d)$ is a term of order $N^d \ell M^{-1} \ll
N^d$, $E_{\psi^N_{t,\mf f}}$ stands for the expectation with respect
to $\psi^N_{t,\mf f} \, d\nu_\alpha^N$, $W_i^*$ is the current in the
$i$-th direction for the adjoint process and $G_{i,j} (\eta)$,
$H_{i,j}(\eta)$ are local functions given by:
\begin{eqnarray*}
\!\!\!\!\!\!\!\!\!\!\!\! && 
W_i^* \;=\; \sum_{y\in\bb Z^d} p^*(y)\, y_i \, \eta(0) [1-\eta(y)]\;, 
\quad G_{i,j} (\eta)\; =\; \sum_{y\in\bb Z^d} p^*(y)\, y_i \, y_j\, 
\eta(0) [1-\eta(y)] \; , \\
\!\!\!\!\!\!\!\!\!\!\!\! &&
\quad H_{i,j}(\eta) \; =\;
\sum_{y\in \bb Z^d} p^*(y) \eta(0) [ 1-\eta (y)] 
\big\{ y_i - \nabla_{0,y} \Gamma_{\mf f_i (\eta^M(x), \cdot)} \big\} \;
\times \\
\!\!\!\!\!\!\!\!\!\!\!\! &&
\qquad\qquad \qquad \qquad \qquad \qquad\qquad \qquad \qquad \qquad
\qquad \qquad 
\;\times \big\{ y_j - \nabla_{0,y} \Gamma_{\mf f_j (\eta^M(x), \cdot)} 
\big\}\;.
\end{eqnarray*}
Here and below, $\nabla_{x,y}$ is the operator defined by
$$
(\nabla_{x,y} f) (\eta) \;=\; f(\sigma^{x,y} \eta) - f(\eta)
$$
and, for a local function $h$, $\Gamma_h$ is the formal sum
$$
\Gamma_h \;=\; \sum_{x\in\bb Z^d} \tau_x h \;.
$$
Since $h$ is a local function, even if the sum of translations is not
defined, the gradient $\nabla_{0,y} \Gamma_h$ makes sense because only
a finite number of terms do not vanish.

We consider separately the sums in (\ref{eq:8}). The goal is to
replace each one by a simpler expression and a remainder denoted
by $o(N^d)$. The remainder $o(N^d)$ stands for an expression
which may depend on time and on the configuration but such that
$$
\lim_{\ell\to\infty}
\limsup_{N\to\infty} N^{-d} \Big\vert \, \int_0^t ds\, o(N^d)  
f^N_s \, d\nu_\alpha^N \, \Big\vert \;=\; 0
$$
for every $t>0$. If the remainder vanishes only after taking the
limit in $\mf f_n$, we denote it by $o(\mf f_n, N^d)$ and we require
$$
\lim_{n\to\infty} \limsup_{\ell\to\infty}
\limsup_{N\to\infty} N^{-d} \Big\vert \, \int_0^t ds\, o(\mf f_n, N^d)  
f^N_s \, d\nu_\alpha^N \, \Big\vert \;=\; 0
$$
for every $t>0$.

We start with the last term of (\ref{eq:8}).  By Corollary \ref{s11},
we may replace the expectation with respect to $\psi^N_{t,\mf f}\,
d\nu_\alpha^N$ with an expectation with respect to
$\nu^N_{\rho(t,\cdot)}$ paying a price of order $N^{d}
\sqrt{|\Lambda_M|/N}$. After this modification, the last line of
(\ref{eq:8}), becomes
$$
- \sum_{x\in\bb T_N^d} (\partial_{t} \lambda) 
(t, x/N) \, \{ \eta (x) - \rho (t, x/N) \}
\;+\; o (N^d)\;.
$$
Since $\partial_{t} \lambda$ is a smooth function, we may further
replace $\eta(x)$ by $\eta^\ell (x)$ paying a price absolutely bounded
by $C \ell^2 N^{d-2}$ for some finite constant $C$.
\smallskip

To estimate the order $N^{d+1}$ term of (\ref{eq:8}), we first
take advantage of the assumption that the solution $\rho(t,u)$ is
constant along the drift direction. 

By paying a price of order $O(\ell^2 N^{d-1})$, we may replace
the current $W_i^*$ by an average $|\Lambda_{\ell}|^{-1}
\sum_{y\in\Lambda_\ell} \tau_y W_i^*$. Here again one should keep
in mind, that the average is in fact carried over a cube of
length slightly smaller than $2\ell +1$ to ensure that all local
functions $\tau_y W_i^*$ have support contained in
$\Lambda_\ell$.

Recall that $q=(q_1, \dots, q_d)$ denotes the drift of particles.
The average of the current $W_i^*$ can be written as
$$
\frac 1{|\Lambda_{\ell'}| }\sum_{z\in\Lambda_\ell'} \tau_z
W_i^*\;=\; q_i^* \Big\{ \eta^{\ell'}(0) - 2 \eta^\ell(0)
\eta^{\ell'}(0) + \eta^\ell(0)^2 \Big\} \;
+\; \frac 1{|\Lambda_{\ell'}| }\sum_{z\in\Lambda_\ell'}
w_i^*(\eta^\ell(0), \tau_z \eta)\; ,
$$
where $q_i^* = -q_i$ and 
$$
w_i^*(\alpha, \eta) \;=\;
- \sum_{y\in\bb Z^d} p^*(y) \, y_i\, [\eta(0)-\alpha]\,
[\eta(y)-\alpha]
\; -\; \alpha \sum_{y\in\bb Z^d} p^*(y) \, y_i\,
[\eta(y)-\eta(0)]\; .
$$
The first term of the current gives no contribution since for any
function $J$,
$$
N \sum_{i=1}^d \sum_{x\in\bb T_N^d} (\partial_{u_i} \lambda) (t, x/N)
q_i J (\eta^{\ell'}(0), \eta^{\ell}(x))\; =\; 0
$$
because $\sum_{1\le i\le d} q_i (\partial_{u_i} \lambda) (t,
u) = \{\rho (t,u) [1-\rho(t,u)]\}^{-1} \sum_{1\le i\le d} q_i
(\partial_{u_i} \rho) (t, u)$ vanishes for all $(t,u)$. The first
term of (\ref{eq:8}) becomes therefore
$$
N \sum_{i=1}^d \sum_{x\in\bb T_N^d} (\partial_{u_i} \lambda) (t, x/N)
\tau_x \Big\{ (A_\ell w_i^*) (\eta^\ell(0), \eta) - 
L_N^* (A_\ell \mf f_i) (\eta^M(0), \eta) \Big\}\; .
$$
To ensure that the function which appears in $A_\ell w_i^*$
has mean zero with respect to the all canonical measures on the cube
$\Lambda_\ell$, we further replace $A_\ell w_i^*$ by $A_\ell^0
w_i^*$, where
$$
(A_\ell^0 w_i^*)(\alpha, \eta) \;=\; (A_\ell w_i^*)(\alpha, \eta)
\; + \; q_i \frac {\alpha (1-\alpha)}{|\Lambda_\ell| - 1}\; \cdot
$$
This replacement is pemitted because $\sum_i q_i \partial_{u_i}
\rho =0$.

Following the nongradient method, we add and subtract $\sum_{1\le
  j\le d} D_{i,j} (\eta^{\varepsilon N}(0))$ $[ \eta^{\varepsilon
  N}(e_j) - \eta^{\varepsilon N}(0)]$. Since the diffusion
coefficient is smooth, this expression is equal to $\sum_{1\le
  j\le d} \{ d_{i,j} (\eta^{\varepsilon N}(e_j)) - d_{i,j}
(\eta^{\varepsilon N}(0)) \} + (\varepsilon N)^{-2}$, where
$d_{i,j}$ stands for the integral of $D_{i,j}$.  In particular,
after a summation by parts, the first line of (\ref{eq:8}), may
be rewritten as
\begin{eqnarray}
\label{eq:12}
\!\!\!\!\!\!\!\!\! &&
N \sum_{i=1}^d \sum_{x\in\bb T_N^d} (\partial_{u_i} \lambda) (t, x/N)
\tau_x  V_{i}^{\varepsilon N, M, \ell} (\eta) \\
\!\!\!\!\!\!\!\!\! && \quad
+\, \sum_{i, j=1}^d \sum_{x\in\bb T_N^d} (\partial^2_{u_i,u_j}
\lambda) (t, x/N) d_{i,j} (\eta^{\varepsilon N}(x)) \; +\;
O(N^{d-1})\; ,
\nonumber
\end{eqnarray}
where
\begin{eqnarray*}
\!\!\!\!\!\!\!\!\!\!\!\!\! &&
V_{i}^{K, M, \ell} (\eta) \;=\;  \\
\!\!\!\!\!\!\!\!\!\!\!\!\! && \quad
(A_\ell^0 w_i^*) (\eta^\ell(0), \eta) \;+\; \sum_{j=1}^d D_{i,j}
(\eta^{K}(0)) [ \eta^{K}(e_j) - \eta^{K}(0)] \;-\;
L_N^* (A_\ell \mf f_i) (\eta^M(0) , \eta)\;.
\end{eqnarray*}

It is not difficult to see that there exists a finite constant $C
(\alpha)$ such that $H(\mu^N \,|\, \nu^N_\alpha) \le C(\alpha) N^d$
for every probability measure $\mu^N$ on $\{0,1\}^{\bb T_N^d}$. In
particular, by the usual two blocks estimate, since $d_{i,j}$ is
Lipschitz continuous, for every $T>0$,
\begin{eqnarray*}
\!\!\!\!\!\!\!\!\!\!\!\!\! &&
\lim_{\ell \to \infty} \limsup_{\varepsilon \to 0}
\limsup_{N\to\infty} \\
\!\!\!\!\!\!\!\!\!\!\!\!\! && \quad
\int_0^T dt\, \int \nu_\alpha^N(d\eta) \, f_t^N (\eta) 
\, N^{-d} \sum_{x\in\bb T_N^d}  \big | \, d_{i,j} (\eta^{\varepsilon N}(x))
- d_{i,j} (\eta^{\ell}(x)) \, \big | \;=\; 0\;.
\end{eqnarray*}
We may therefore replace in the second line of (\ref{eq:12}) the
average of particles over a small macroscopic cube by the average over
a large microscopic cube, i.e., replace $\eta^{\varepsilon N}(x)$ by 
$\eta^{\ell}(x)$.

On the other hand, the usual nongradient techniques, based on
integration by parts formula, allows the replacement in (\ref{eq:12})
of $D_{i,j} (\eta^{\varepsilon N}(0)) [ \eta^{\varepsilon N}(e_j) -
\eta^{\varepsilon N}(0)]$ by $D_{i,j} (\eta^{\ell}(0)) [
\eta^{\ell'}(e_j) - \eta^{\ell'}(0)]$.  Here $\ell' = \ell -1$ for the
previous function to depend only on the sites in $\Lambda_\ell$. To
keep notation simple, we will denote this expression by $D_{i,j}
(\eta^{\ell}(0)) [ \eta^{\ell}(e_j) - \eta^{\ell}(0)]$. We refer to
Chap. 7 of \cite{kl} for a proof of this replacement.

In subsection \ref{sec4.2} we prove that we may replace $L_N^* (A_\ell
\mf f_i) (\eta^M(0) , \eta)$ by $L^*_{\Lambda_\ell} (A_\ell \mf f_i)$
$(\eta^\ell(0) , \eta)$. Here $L^*_{\Lambda_\ell}$ stands for the
restriction of the generator $L_N^*$ to the cube $\Lambda_\ell$.  This
means that we suppress all jumps from $\Lambda_\ell$ to
$\Lambda_\ell^c$ and all jumps from $\Lambda_\ell^c$ to
$\Lambda_\ell$. In particular, this generator leaves $\eta^\ell(0)$
invariant and it is acting in fact only on the second coordinate.
This replacement is one of the main technical point of the article. It
is here that the special form of $\psi^N_{t,\mf f}$ plays an important
role, that we need the spatial averages and the particular size of $M$
and $\ell$ presented in (\ref{eq:16}).

Up to this point, we transformed the first line of (\ref{eq:8}) in
\begin{eqnarray}
\label{eq:10}
\!\!\!\!\!\!\!\!\! &&
N \sum_{i=1}^d \sum_{x\in\bb T_N^d} (\partial_{u_i} \lambda) (t, x/N)
\tau_x  V_{i}^{\ell} (\eta) \\
\!\!\!\!\!\!\!\!\! && \quad
+\, \sum_{i, j=1}^d \sum_{x\in\bb T_N^d} (\partial^2_{u_i,u_j}
\lambda) (t, x/N) d_{i,j} (\eta^{\ell}(x)) \; +\;
o(N^{d})\; ,
\nonumber
\end{eqnarray}
where
$$
V_{i}^{\ell} (\eta) \;=\;  
(A_\ell^0 w_i^*) (\eta^\ell(0), \eta) \;+\; \sum_{j=1}^d D_{i,j}
(\eta^{\ell}(0)) [ \eta^{\ell}(e_j) - \eta^{\ell}(0)] \;-\;
L^*_{\Lambda_\ell} (A_\ell \mf f_i) (\eta^\ell(0) , \eta) \;.
$$
By the nongradient method, the first line can be shown to be of order
$o(\mf f, N^d)$. Details are given in subsection \ref{sec4.3}.
\smallskip

It remains to consider the second and third line of (\ref{eq:8}).
By the one block estimate the second line of (\ref{eq:8}) is equal to 
$$
(1/2) \sum_{i, j =1}^d \sum_{x\in\bb T_N^d} 
(\partial^2_{u_i, u_j} \lambda) (t, x/N) \, \sigma_{i,j} \tau_x
F(\eta^\ell(0)) \; +\; o(N^{d})\;,
$$
where $\sigma_{i,j}$ is the symmetric matrix defined just after
(\ref{eq:4}) and $F(a) = a(1-a)$.  For $1\le i,j\le d$, let
\begin{equation}
\label{eq:31}
J_{i,j} (\beta) \;=\; 2 \beta (1-\beta) \Big\{ 
D_{i,j} (\beta) - \beta \sigma_{i,j} \Big\} \; .
\end{equation}
We prove in subsection \ref{sec4.1} that the third line of
(\ref{eq:8}) is equal to
\begin{eqnarray*}
\!\!\!\!\!\!\!\!\!\!\!\!\! &&
(1/2) \sum_{i,j=1}^d \sum_{x\in\bb T_N^d} (\partial_{u_i}
\lambda) (t, x/N) (\partial_{u_j} \lambda) (t, x/N) 
\sigma_{i,j} F(\eta^\ell (x)) \\
\!\!\!\!\!\!\!\!\!\!\!\!\! && \quad +\;
(1/2) \sum_{i,j=1}^d \sum_{x\in\bb T_N^d} (\partial_{u_i}
\lambda) (t, x/N) (\partial_{u_j} \lambda) (t, x/N) 
J_{i,j} (\eta^\ell (x)) \; +\; o(\mf f_n, N^d)\; .
\end{eqnarray*}

In conclusion, we proved that (\ref{eq:8}) is equal to
\begin{eqnarray}
\label{eq:13}
\!\!\!\!\!\!\!\!\!\!\!\!\! && 
\sum_{m=1}^2 \sum_{i,j=1}^d \sum_{x\in\bb T_N^d} G^m_{i,j}(t,
x/N) H^m_{i,j} (\eta^\ell (x)) \\
\!\!\!\!\!\!\!\!\!\!\!\!\! && \quad \nonumber
-\; \sum_{x\in\bb T_N^d} (\partial_{t} \lambda) 
(t, x/N) \, \{ \eta^\ell (x) - \rho (t, x/N) \}
\;+\; o (\mf f_n, N^d)\;.
\end{eqnarray}
where
\begin{eqnarray*}
\!\!\!\!\!\!\!\!\!\!\!\!\! &&
G^1_{i,j}(t,u) = (\partial^2_{u_i,u_j} \lambda)(t,u)\;, \quad
H^1_{i,j} (\beta) \; =\; d_{i,j}(\beta) \; +\; (1/2) \sigma_{i,j}
F(\beta)\;, \\
\!\!\!\!\!\!\!\!\!\!\!\!\! && \quad
G^2_{i,j}(t,u) = (\partial_{u_i} \lambda)(t,u) (\partial_{u_j}
\lambda)(t,u) \;, \quad
H^2_{i,j} (\beta) \; =\; (1/2) \{ J_{i,j}(\beta) \; +\; \sigma_{i,j}
F(\beta)\} \;.
\end{eqnarray*}

An integration by parts shows that
\begin{equation*}
\sum_{m=1}^2 \sum_{i, j=1}^d \int_{\bb T^d} du\, 
G^m_{i,j}(t,u) H^m_{i,j} (\rho(t,u)) \;=\; 0 \;.
\end{equation*}
In particular, in formula (\ref{eq:13}), we may replace the
terms $H^m_{i,j} (\eta^\ell (x))$ by $H^m_{i,j}(\eta^\ell
(x)) - H^m_{i,j} (\rho(t,x/N))$ paying a price of order $o(N^d)$.
A further elementary computation gives that
\begin{equation*}
\sum_{m=1}^2 \sum_{i,j=1}^d G^m_{i,j}(t, u) 
(H^m_{i,j})' (\rho(t,u)) \;=\;  (\partial_{t} \lambda) 
(t, u) 
\end{equation*}
for every $t$ and $u$, where $(H^m_{i,j})'$ stands for the
derivative of $H^m_{i,j}$. Therefore, (\ref{eq:13}) becomes
\begin{equation*}
\sum_{m=1}^2 \sum_{i,j=1}^d \sum_{x\in\bb T_N^d} 
G_{i,j}^m (t, x/N) \, B_{i,j}^m  (\eta^{\ell}(x),\rho(t, x/N))
\;+\; o (\mf f_n, N^d)\;,
\end{equation*}
where
$$
B_{i,j}^m (a,b) \;=\; H^m_{i,j}(a) \;-\; 
H^m_{i,j} (b) \;-\; (H^m_{i,j})' (b) \, [ a-b ] \;.
$$
At this point we may repeat the standard arguments of the relative
entropy method do conclude. We refer to Chap. 6 of \cite{kl} for
details.

\section{Hilbert space of variances}
\label{sec6}

We prove in this section the existence of functions $\mf f_1,
\dots \mf f_d$ in $\mf F$ which approximate the current in the
Hilbert space of variances. We rely on recent results based on
general duality presented in \cite{lov1}, \cite{lov2}.

For $0\le \alpha\le 1$, denote by $\mc G_\alpha$ the space of
cylinder functions $g$ such that $E_{\nu_{\alpha}} [g] =
\partial_\alpha E_{\nu_{\alpha}} [g] =0$:
$$
\tilde g (\alpha) \; =\; E_{\nu_{\alpha}} [g] \; =\; 0\quad
\text{and}\quad \tilde g' (\alpha)\; =\; \frac{d}{d\beta} 
E_{\nu_{\beta}} [g] \Big\vert_{\beta=\alpha}\; =\;  0\; .
$$
For each function $g$ in $\mc G _\alpha$ we define $\V
g\V_\alpha$ by
\begin{equation}
\label{eq:30}
\V g\V_\alpha^2 \;=\; |g|_\alpha^2 \; +\; \Vert g\Vert_{-1, \alpha}^2\;,
\end{equation}
where
$$
\aligned
& |g|_\alpha^2  \; =\; \sup _{a \in \bb R^d} \Big\{ 2 \sum_{i=1}^d a_i \, 
\sum_{x\in\bb Z^d} x_i  < g \,;\, \eta (x)>_\alpha -
\frac{\chi (\alpha)}{2} a \cdot \sigma a  \Big\}\; , \\
& \quad \Vert g\Vert_{-1, \alpha}^2 \; =\; \sup_{h \in \mc G_\alpha}  \Big\{
2 \ll g,h \gg_\alpha  - \ll h , (-L^s) h \gg_\alpha \Big\}\; .
\endaligned
$$
In this formula, $\chi(\alpha) = \alpha(1-\alpha)$, $a\cdot b$
stands for the inner product in $\bb R^d$ and $\ll \cdot, \cdot
\gg_\alpha$ for the inner product in $\mc G_\alpha$ given by
$$
\ll g,h \gg_\alpha \; =\; \sum_{x\in\bb Z^d} 
< g \,;\, \tau_x h>_\alpha\;,
$$
where $< f_1 \,;\, f_2>_\alpha$ denotes the covariance of $f_1$,
$f_2$ with respect to $\nu_\alpha$. Notice that in the sums which
appear in the formulas above, all but a finite number of terms vanish
because $\nu_\alpha$ is a product measure. Theorem \ref{s7} is the
main result of this section.

\begin{theorem}
\label{s7}
There exist a smooth matrix-valued function $D(\alpha) = \{D_{i,j}
(\alpha),\, 1\le i,j\le d\}$ and a sequence of functions $\{\mf
f_{i,n},\, n\ge 1\}$ in $\mf F$, $1\le i\le d$, such that
$$
\lim_{n\to\infty} \sup_{\alpha\in [0,1]}
\V w_i^*(\alpha, \eta) + \sum_{j=1}^d D_{i,j}(\alpha)
[\eta(e_j) - \eta(0)] - L^* \mf f_{i,n} (\alpha, \eta) 
\V_{\alpha} \;=\;0
$$
for $1\le i\le d$. Moreover, for any vector $v$ in $\bb R^d$,
\begin{eqnarray}
\label{eq:26}
\!\!\!\!\!\!\!\!\!\!\!\!\! &&
\lim_{n\to\infty} \sum_{x\in\bb Z^d} 
< \sum_{j=1}^d v_j \mf f_{j ,n} (\alpha, \eta), (-L^s) 
\tau_x \sum_{j=1}^d v_j \mf f_{j ,n} (\alpha, \eta) >_\alpha \\
\!\!\!\!\!\!\!\!\!\!\!\!\! && \nonumber \quad 
\;=\; \chi(\alpha) \, v \cdot \{ D(\alpha) - \alpha \sigma\} v
\end{eqnarray}
uniformly in $\alpha$.
\end{theorem}

This result is a slight generalization of Corollary 10.1 and Lemma
10.4 in \cite{loy1}, proved in \cite{loy2} using results presented in
\cite{ly}. We have the advantage here to obtain uniformity up to the
boundary. In sake of completeness, we present a simpler proof relying
on the generalized duality developed in \cite{lov1}, \cite{lov2}.

To keep notatiom simple, we prove Theorem \ref{s7} for the
current $w_i$ obtained from $w_i^*$ by replacing $p^*(\cdot)$ by
$p(\cdot)$ and for the generator $L$ in place of $L^*$.

\noindent{\bf Duality.} For each $n\ge 0$, denote by $\mc E_{n}$ the
subsets of $\mathbb Z^d$ with $n$ points and let $\mc E = \cup_{n\ge
  0} \mc E_{n}$ be the class of finite subsets of $\mathbb Z^d$. For
each $A$ in $\mc E$, let $\Psi_A$ be the local function
$$
\Psi_A \;=\;  \prod_{x\in A} \frac{ \eta (x) -
  \alpha}{\sqrt{\chi(\alpha)}}\; .
$$
By convention, $\Psi_\phi = 1$. It is easy to check that
$\{\Psi_A,\, A\in \mc E \}$ is an orthonormal basis of $L^2
(\nu_\alpha)$.  For each $n\ge 0$, denote by $\mc D_n$ the subspace of
$L^2(\nu_\alpha)$ generated by $\{\Psi_A, \, A\in \mc E_{n}\}$, so
that $L^2(\nu_\alpha)= \oplus_{n\ge 0} \mc D_n$. Functions in $\mc
D_n$ are said to have degree $n$.

Consider a local function $f$. Since $\{\Psi_A :\, A \in \mc E\}$ is a
basis of $L^2(\nu_\alpha)$, we may write
$$
f=\sum_{n\ge 0} \sum_{A\in \mc E_{n}} \mathfrak f(\alpha, A) \Psi_A \; .
$$
Note that the coefficients ${\mf f}(\alpha, A)$ depend not
only on $f$ but also on the density $\alpha$. Since $f$ is a
local function, $\mathfrak f \colon \mc E \to \bb R$ is a
function of finite support.

Fix a local function $f$ and denote by $\mf f(\alpha, A)$ its
Fourier coefficients. $f$ has zero mean with respect to
$\nu_\alpha$ if and only if $\mf f(\alpha, \phi)=0$. It belongs
to $\mc G_\alpha$ if and only if $\mf f(\alpha, \phi)=0$ and the
degree one part is such that
$$
\sum_{z\in\bb Z^d} \mf f(\alpha, \{z\}) \;=\; 0\; .
$$
In this case, we may rewrite the degree one piece as
$$
\sqrt{\chi(\alpha)} \sum_{z\in\bb Z^d} \mf f(\alpha, \{z\}) 
[\eta(z) - \eta(0)]\;.
$$
In conclusion, all functions $f$ in $\mc G_\alpha$ may be written
as
$$
\sqrt{\chi(\alpha)} \sum_{z\in\bb Z^d} \mf f(\alpha, \{z\})  
[\eta(z) - \eta(0)]\;+\;
\sum_{n\ge 2} \sum_{A\in \mc E_{n}} \mathfrak f(\alpha, A)
\Psi_A\; .
$$
For $n\ge 0$, denote by $\pi_n$ the projection on $\mc D_n$ so that
$f = \sum_{n\ge 1} \pi_n f$ for $f$ in $\mc G_\alpha$. In the formula
above, the first term corresponds to $\pi_1 f$, the piece of $f$ which
has degree one, and the second term corresponds to $(I-\pi_1) f$, the
piece of degree greater or equal to $2$.

It is clear that a local function of type $h- \tau_x h$ belongs to the
kernel of the inner product $\ll \cdot, \cdot\gg_\alpha$ defined
above. This is the case of $\eta(z) - \eta(0)$ so that $\Vert
f\Vert_{-1,\alpha} = \Vert (I-\pi_1) f \Vert_{-1,\alpha}$.  In
contrast, any function $h$ of degree greater or equal to $2$ is such
that
$$
\sum_{x\in\bb Z^d} x_i  < h \,;\, \eta (x)>_\alpha \; =\; 0
$$
for all $i$ so that $|h|_{\alpha}=0$. Therefore, $|f|_\alpha =
|\pi_1 f|_\alpha$ and 
$$
\V f\V^2_\alpha \;=\; |\pi_1 f|_\alpha^2 \;+\;
\Vert (I-\pi_1) f \Vert_{-1,\alpha}^2
$$
for every local function $f$ in $\mc G_\alpha$.

\noindent{\bf The generator on the Fourier coefficient.} Let $\mc E_*$ be
the class of all finite subsets of $\bb Z^d_*= \bb Z^d \backslash
\{0\}$ and let $\mc E_{*,n}$ be the class of all subsets of $\bb
Z^d_*$ with $n$ points. For a local function $f$ in $\mc G_\alpha$,
define $\mf T {f} : [0,1]\times \mc E_* \to \bb R$ by
$$
(\mf T {f}) (\alpha, A) \;=\;
\sum_{z\in\bb Z^d} {\mf f}(\alpha, [A\cup \{0\}]+z) \;,
$$
where ${\mf f}(\alpha, B)$ stands for the Fourier coefficients of
$f$.  In this context, a function $f(\alpha, \eta)$ belongs to $\mc
G_\alpha$ if and only if $\mf f(\alpha, \phi) = (\mf T f)(\alpha,
\phi)=0$.  It has been in proved in \cite{lov2} that for every
zero-mean local functions $f$, $g$
\begin{equation}
\label{eq:23}  
\ll f,g \gg_\alpha \; = \; <(\mf T {f}) , (\mf T {g})> \;=\;
\sum_{n\ge 0} \frac 1{n+1}
\sum_{A\in \mc E_{*,n}} (\mf T {f}) (\alpha, A) \, (\mf T {g})
(\alpha,  A)\; .
\end{equation}
For functions in $\mc G_\alpha$, this sum starts from $1$ because
$(\mf T f)(\alpha, \phi) = (\mf T g)(\alpha,\phi)=0$.

Observe that not every function ${\mf f} : [0,1]\times \mc E_* \to
\bb R$ is the image by $\mf T$ of some local function $f$ since
\begin{equation}
\label{eq:22}
(\mf T f) (\alpha, A) \; =\; (\mf T f) (\alpha, S_z A)
\end{equation}
for all $z$ in $A$. Here, $S_z A$ is the set defined by
\begin{equation*}
S_z A \;=\; \left\{
\begin{array}{ll}
A-z & \text{ if $z\not \in A$,} \\
(A-z)_{0,-z} & \text{ if $z\in A$}\; .
\end{array}
\right.
\end{equation*}

Let ${\mf f}_* : [0,1]\times \mc E_* \to \bb R$ be a finitely
supported function satisfying (\ref{eq:22}).  Define ${\mf f} :
[0,1]\times \mc E \to \bb R$ by
\begin{equation}
\label{eq:24}
{\mf f} (\alpha, B) \;=\; \left\{
\begin{array}{cl}
|B|^{-1} {\mf f}_* (\alpha, B \setminus \{0\}) & \text{if $B\ni 0$}\; , \\
0 & \text{otherwise}\; .
\end{array}
\right.
\end{equation}
An elementary computations shows that $\mf T f(\alpha, \eta) = {\mf
  f}_*$, if $f(\alpha, \eta)$ is the local function whose Fourier
coefficients are $\mf f(\alpha, A)$. Notice that $f(\alpha, \eta)$
belongs to $\mc G_\alpha$ if $\mf f_* (\alpha, \phi)=0$.

For any local function $f$, $\mf T (L f) = {\mf L}_\alpha \mf T {f}$,
provided
\begin{equation*}
{\mf L}_\alpha = {\mf L}_s + (1-2\alpha) {\mf L}_d
+ \sqrt{\chi(\alpha)} \{ {\mf L}_+ + {\mf L}_- \}
\end{equation*}
and, for $A\in \mathcal E_*$, $\mf v : \mc E_* \to\bb R$ a finitely
supported function,
\begin{eqnarray*}
\!\!\!\!\!\!\!\!\!\!\!\!\!\!\!\!\!\!\!\! &&
(\mf L_{s} {\mf v}) \, (B) \;=\; (1/2)
\sum_{x, y\in \mathbb Z^d_*}  s(y-x)
[ {\mf v} (B_{x,y}) - {\mf v}(B)] \; +\; 
\sum_{y\not\in B} s(y) [ {\mf v} (S_y B) - {\mf v}(B) ]\;, \\
\!\!\!\!\!\!\!\!\!\!\!\!\!\!\!\!\!\!\!\! && \quad
({\mf L}_d {\mf v}) (A) \;=\; \sum_{
\substack{x\in A, y\not\in A \\ x,y\neq 0}} a(y-x)
\{ {\mf v}(A_{x,y}) - {\mf v}(A))
+ \sum_{\substack{ y\not\in A \\ y\neq 0}} a(y)
\{ {\mf v}( S_y A) - {\mf v}(A)\} \; , \\
\!\!\!\!\!\!\!\!\!\!\!\!\!\!\!\!\!\!\!\! && 
\quad \quad ({\mf L}_+ {\mf v}) (A) \;=\;
2 \sum_{x\in A, y\in A} a(y-x) \, {\mf v} (A\backslash \{y\}) \\
\!\!\!\!\!\!\!\!\!\!\!\!\!\!\!\!\!\!\!\!
&&\qquad\qquad\qquad\qquad
 +\;  2 \sum_{x\in A} a(x) \{ {\mf v}(A\backslash \{x\})
- {\mf v}(S_x [A\backslash \{x\}]) \} \; , \\
\!\!\!\!\!\!\!\!\!\!\!\!\!\!\!\!\!\!\!\!
&& \qquad\quad ({\mf L}_- {\mf v}) (A) \;=\; 
2 \sum_{\substack {x\not\in A, y\not\in A\\ x,y\neq 0}} a(y-x)
\, {\mf v}(A\cup \{y\}) \; .
\end{eqnarray*}
In this formula, $A_{x,y}$ is the set defined by
\begin{equation*}
A_{x,y} \; =\;
\left\{
\begin{array}{ll}
(A \backslash \{x\}) \cup \{y\} & \text{if $x\in A$, $y \not \in A$,} \\
(A  \backslash \{y\}) \cup \{x\} & \text{if $y\in A$, $x \not \in A$,} \\
A & \text{otherwise}\; ;.
\end{array}
\right.
\end{equation*}

\noindent{\bf Hilbert spaces.} For two local functions
$f$, $g$, let
$$
\ll f , g\gg_{\alpha, 1} \;=\; \ll f, (-L^s) g\gg_\alpha
$$
and let $H_1(\alpha)$ be the Hilbert space generated by local
functions $f$ and the inner product $\ll \cdot , \cdot \gg_{\alpha,
  1}$.  Denote by $\ll \cdot, \cdot\gg_1$ the scalar product on $\mc
E_{*}$ defined by
\begin{equation*}
\ll \mf f, \mf g \gg_{1} \; =\; 
\sum_{n\ge 0} \frac 1{n+1} \sum_{A\in \mc E_{*,n}}
{\mf f} (\alpha, A)  (- {\mf L}_{s} {\mf g} ) (\alpha, A)
\end{equation*}
and by $\mf H_1$ the Hilbert space generated by the finite supported
functions endowed with the previous scalar product.  From the previous
definitions, for every local function $f$, $g$,
$$
\ll f,g \gg_{1,\alpha} \;=\; \ll \mf T f , \mf T g \gg_{1}
$$

To introduce the dual Hilbert spaces of $H_1$, $\mf H_1$, for a local
function $f$, consider the semi-norm $\Vert \,\cdot\,\Vert_{-1}$ given
by
$$
\Vert f\Vert_{-1, \alpha}^2 \;=\; \sup_{g} \Big\{ 2 \ll f, g \gg_\alpha 
- \ll g, g \gg_{1,\alpha} \Big\}\; ,
$$
where the supremum is carried over all local functions $g$.  Denote
by $H_{-1}$ the Hilbert space generated by the local functions and the
semi-norm $\Vert \cdot\Vert_{-1}$.  In the same way, for a finitely
supported function $\mf f : \mc E_{*} \to\bb R$, let
$$
\Vert \,\mf f\, \Vert_{-1}^2 \;=\;
\sup_{\mf g} \Big\{ 2 <\mf f, \mf g> - <\mf g, \mf g>_1 \Big\}\; ,
$$
where the supremum is carried over all finitely supported functions
$\mf g : \mc E_{*} \to \bb R$ and $<\cdot , \cdot>$ is the inner
product on $L^2(\mc E^*)$ defined in (\ref{eq:23}). Denote by $\mf
H_{-1}$ the Hilbert space induced by the finitely supported functions
$\mf f : \mc E_{*} \to\bb R$ and the semi-norm $\Vert \cdot
\Vert_{-1}$. By the identities for the $L^2$ and the $H_1$ norms, we
obtain that
\begin{equation}
\label{eq:28}
\Vert f(\alpha, \eta) \Vert_{-1, \alpha}^2 \;=\;
\Vert (\mf T f) (\alpha, \cdot) \Vert_{-1}^2
\end{equation}

\noindent{\bf The currents.} Recall the definition of the current
$w_i(\alpha, \eta)$ given in section \ref{sec5}.  $w_i$ is
obtained from $w_i^*$ by replacing $p^*(\cdot)$ by $p(\cdot)$ and
can be expressed as
$$
w_i\;=\; - \alpha (1-\alpha) \sum_{y\in\bb Z^d} p(y)\, y_i\,
\Psi_{0,y} \;-\; \alpha \sum_{y\in\bb Z^d} p(y)\, y_i\, \{ \eta (y) -
\eta(0)\}\; .
$$
Denote the first piece, which has degree $2$, by $\alpha (1-\alpha)
w_i^0$. On the other hand, since $\eta(e_k) - \eta(0) = \eta(e_k+x) -
\eta (x)$ for the norm $|\cdot|_\alpha$ for any $x$, the piece which
has degree one is equal to $\alpha \sum_{y\in\bb Z^d} \sum_{1\le j\le
  d} p(y)\, y_i\, y_j \{ \eta (e_j) - \eta(0)\}$ so that
$$
w_i\;=\; \alpha (1-\alpha) w_i^0 \;-\; \alpha \sum_{j=1}^d
\sigma_{i,j} \, [\eta (e_j) - \eta(0)] \; . 
$$
Let $\mf w_i = \mf T w_i^0$. An straightforward computation gives
that
$$
\mf w_i (\alpha, \{z\})  \;=\; - 2 \, z_i \, a(z)
$$
for $z\not = 0$ and $\mf w_i (\alpha, A)=0$ otherwise. Notice that
$\mf w_i$ does not depend on $\alpha$.  We have now all elements to
prove the main result of this section.

\noindent{\bf Proof of Theorem \ref{s7}.} Fix $1\le i\le d$. By
Theorem 4.1 in \cite{lov2}, $\mf w_i$ belongs to $\mf H_{-1}$ because
$\mf w_i (\alpha, \phi)=0$ and we are in $d\ge 3$.

It has been proved in Lemma 4.3 of \cite{lov2} that for each $\lambda
>0$ there exists a solution $\mf f_{i,\lambda} (\alpha, A)$ of the
resolvent equation
$$
\lambda \mf f_{i,\lambda} \;-\; \mf L_\alpha \mf f_{i, \lambda} \;=\; 
\mf w_i
$$
satisfying (\ref{eq:22}) and such that $\mf f_{i,\lambda}(\alpha,
\phi)=0$. 

By Theorem 4.4 in \cite{lov2}, for any $k\ge 1$, there exists a finite
constant $C_k$ independent of $\alpha$ and $\lambda$ such that
\begin{equation}
\label{eq:11}
\lambda \sum_{n\ge 0} (1+n)^k <\pi_n \mf f_{i,\lambda} ,
\pi_n \mf f_{i,\lambda}>_n \;+\;
\sum_{n\ge 0} (1+n)^k < \pi_n \mf f_{i,\lambda} 
, (- \mf L_s) \pi_n \mf f_{i,\lambda}>_n \;\le\; C_k
\end{equation}
for every $\lambda>0$ and $\alpha$ in $[0,1]$. In this formula,
$\pi_n$ stands for the projection on $\mc E_{*,n}$: $(\pi_n \mf
f)(\alpha, A) = \mf f (\alpha, A) \mb 1\{A\in \mc E_{*,n}\}$, and
$<\cdot , \cdot>_n$ for the inner product in $\mc E_{*,n}$ with
respect to the counting measure:
$$
< \mf f , \mf g >_n \;=\; \sum_{A\in \mc E_{*,n}} \mf f (\alpha, A)
\mf g (\alpha, A)\;.
$$
The estimate is uniform in $\alpha$ because the current $\mf w_i$
does not depend on $\alpha$.

By section 6 of \cite{lov2}, for each $z$ in $\bb Z^d$, $\mf
f_{i,\lambda} (\cdot, \{z\})$ is a smooth function in $[0,1]$ and
there exists a subsequence $\lambda_k \downarrow 0$ such that $\mf
f_{i,\lambda_k} (\alpha, \{z\})$ converges uniformly, as well as its
derivatives, to some smooth function $\mf f_i (\alpha, \{z\})$.

By the proof of Lemma 2.8 of \cite{lov3}, taking a further
subsequence, we may assume that $- (\mf L_\alpha \mf f_{i,\lambda_k})
(\alpha, \cdot)$ converges weakly to $\mf w_i$ in $\mf H_{-1}$ for a
countable dense subset of densities $\alpha$ in $[0,1]$. By Lemma
\ref{s12} below, $- (\mf L_\alpha \mf f_{i,\lambda_k}) (\alpha,
\cdot)$ converges weakly to $\mf w_i$ in $\mf H_{-1}$ for all $\alpha$
in $[0,1]$.

Our goal is to replace the sequence $\mf f_{i,\lambda_k}$ by a
sequence $\mf h_{i,n}$ of finite supported functions with all the
above properties of $\mf f_{i,\lambda_k}$ and for which $- (\mf
L_\alpha \mf h_{i,n}) (\alpha, \cdot)$ converges strongly to $\mf w_i$
in $\mf H_{-1}$, uniformly in $\alpha$.

For each $\alpha$ fixed, we may take convex combinations of the
functions $\mf f_{i,\lambda_k}$ to obtain a new sequence $\mf g_{i,n}$
such that $- \mf L_\alpha \mf g_{i,n}$ converges strongly to $\mf w_i$
in $\mf H_{-1}$. Lemma \ref{s12} below shows that the procedure can
be made uniform in $\alpha$.  Indeed, fix $\varepsilon >0$ and a
finite set $\{\alpha_j, 1\le j\le m\}$ in $[0,1]$. The standard
procedure to derive a strong converging sequence from a weak, bounded
converging sequence shows that there exist $M\ge 1$ and a probability
$(\theta_1, \dots, \theta_M)$ in $\{1, \dots, M\}$, such that
$$
\max_{1\le j\le m} \Vert \mf L_{\alpha_j} \mf g (\alpha_j, \cdot)
+ \mf w_i \Vert_{-1} \;\le\; \varepsilon\;,
$$
where
$$
\mf g (\alpha_j, \cdot) \;=\; \sum_{l=1}^M \theta_l \, \mf
f_{i,\lambda_l} (\alpha_j, \cdot) \; .
$$
Notice that we are taking the same convex combination for all
densities $\alpha_j$.  If $m$ is equal to $\delta^{-1}$, given by
Lemma \ref{s12} below, and $\alpha_j = j \delta$, by Lemma \ref{s12},
$$
\sup_{\alpha\in [0,1]} \Vert \mf L_{\alpha} \mf g (\alpha, \cdot)
+ \mf w_i \Vert_{-1} \;\le\; 2 \varepsilon\;,
$$
where $\mf g (\alpha, \cdot)$ is obtained from $\mf f_{i,\lambda}
(\alpha, \cdot)$ through the same convex combination. We have thus
constructed a convex combination which guarantees the strong
convergence in $\mf H_{-1}$ for all values of $\alpha$. That is, there
exists a sequence $\mf g_{i,n} (\alpha, \cdot)$ such that
\begin{itemize}
\item For each $n\ge 1$, and each $z$ in $\bb Z^d$, $\mf g_{i,n}
  (\cdot, \{z\})$ is a smooth function of $\alpha$ which converges
  uniformly, as well as all its derivatives, to some smooth function
  $\mf f_i (\alpha, \{z\})$.
  
\item Each $\mf g_{i,n}$ satisfies (\ref{eq:22}) and $\mf g_{i,n}
  (\alpha, \phi)=0$ because the functions $f_{i,\lambda_k}$ have this
  property.

\item The sequence converges uniformly to $\mf w_i$ in $\mf H_{-1}$:
$$
\lim_{n\to\infty} \sup_{\alpha\in [0,1]} 
\Vert \mf L_{\alpha} \mf g_{i,n} (\alpha, \cdot) + \mf w_i \Vert_{-1}
\;=\; 0\;.
$$
\end{itemize}

It remains to replace the functions $\mf g_{i,n}$ by finite supported
functions.  Fix two integer $m$, $\ell$ and let $\mf h_{i,n} (\alpha,
A) = \mf g_{i,n} (\alpha, A) \mb 1\{ |A|\le m, A < \Lambda_\ell\}$.
The integers $m$, $\ell$, which depend on $n$ and increase to infinity
with $n$, will be chosen later. Here, $A < \Lambda_\ell$ if there
exists $z$ in $A$ such that $S_z A \subset \Lambda_\ell$. In this way,
$\mf h_{i,n}$ satisfies (\ref{eq:22}).

The sequence $\mf h_{i,n}$ just defined has the first two properties
of the sequence $\mf g_{i,n}$ enumerated above because $m$ and
$\ell$ increase to infinity as $n\uparrow\infty$. To prove the third
one, recall from the computations performed after (4.12) in \cite{loy2}
that
\begin{eqnarray*}
\!\!\!\!\!\!\!\!\!\!\!\!\! &&
\Vert \mf L_{\alpha} \mf g_{i,n} (\alpha, \cdot) -
\mf L_{\alpha} \mf h_{i,n} (\alpha, \cdot) \Vert_{-1}^2 \\
\!\!\!\!\!\!\!\!\!\!\!\!\! && \quad
\le\; C_0 \sum_{k=1}^{m+1} (1+k) \, \Vert \pi_k \mf g_{i,n} (\alpha,
\cdot) -  \pi_k \mf h_{i,n} (\alpha, \cdot) \Vert_{0,k}^2 \\
\!\!\!\!\!\!\!\!\!\!\!\!\! && \qquad
+\; \sum_{k\ge m} (1+k)^2 < \pi_k \mf g_{i,n} (\alpha,
\cdot) , (\mf L_s \pi_k \mf g_{i,n}) (\alpha, \cdot)>_k
\end{eqnarray*}
for some finite constant $C_0$ independent of $\alpha$. Here, $\Vert
\cdot \Vert_{0,k}$ stands for the norm associated to the scalar
product $< \cdot, \cdot >_k$ defined above.  By (\ref{eq:11}), the
second term on right hand side can be made uniformly small in $\alpha$
by choosing $m$ large enough because each function $\mf g_{i,n}$ is
obtained as convex combinations of the solution of the resolvent
equation. For a fixed finite set $\alpha_1, \dots, \alpha_r$, we may
turn the first term as small as one wishes for $\{\alpha_i \, 1\le
i\le r\}$ by taking $\ell$ large enough. By Lemma \ref{s12} below, we
may turn the estimate uniform in $\alpha$ because the functions $\mf
g_{i,n}$ are convex combinations of the solution of the resolvent
equation.

For each fixed $n$, the functions $\mf h_{i,n}( \alpha, \cdot)$ has a
uniform support. Since $\mf h_{i,n}$ satisfies (\ref{eq:22}) and $\mf
h_{i,n} (\alpha, \phi)=0$, the local functions $f_{i,n} (\alpha,
\eta)$ obtained from $\mf h_{i,n}$ through (\ref{eq:24}) are in $\mf
F$.

We claim that the sequence $- \chi(\alpha) f_{i,n} (\alpha, \eta)$ has
the properties required in the statement of the theorem. In view of
the decomposition of the current $w_i$, by (\ref{eq:30}),
\begin{eqnarray}
\label{eq:25}
\!\!\!\!\!\!\!\!\!\!\!\! &&
\Big |\! \Big| \! \Big| w_i(\alpha, \eta) + 
\sum_{j=1}^d D_{i,j}(\alpha) [\eta(e_j) - \eta(0)] + \chi(\alpha) L
f_{i,n} (\alpha, \eta) \Big |\! \Big| \! \Big|_\alpha^2 \\
\!\!\!\!\!\!\!\!\!\!\!\! && \nonumber \quad =\;
\Big \Vert \chi(\alpha) w_i^0 + \chi(\alpha) (I-\pi_1)  L
f_{i,n} (\alpha, \eta) \Big \Vert_{-1, \alpha}^2 \\
\!\!\!\!\!\!\!\!\!\!\!\! && \nonumber \qquad +\;
\Big |  \sum_{j=1}^d \{ D_{i,j}(\alpha) - \alpha \sigma_{i,j}\}
[\eta(e_j) - \eta(0)] + \chi(\alpha) \pi_1  L
f_{i,n} (\alpha, \eta) \Big |_\alpha^2\; .
\end{eqnarray}
Since functions of degree $1$ are in the kernel of the scalar product
$\ll \cdot, \cdot \gg_\alpha$, we may replace $(I-\pi_1) L f_{i,n}$ by
$ L f_{i,n}$ on the first term on the right hand side. On the other
hand, by definition of $\mf T$, by identity (\ref{eq:28}) and since
$\mf T w_i^0= \mf w_i$, the first term on the right hand side of
(\ref{eq:25}) is equal to
$$
\chi(\alpha)^2 \, \big \Vert \mf w_i + \mf L_\alpha
\mf h_{i,n} (\alpha, \cdot) \big \Vert_{-1}^2\;.
$$
This expression vanishes, as $n\uparrow\infty$, uniformly in $\alpha$,
by construction of the sequence $\mf h_{i,n}$.

On the other hand, an elementary computation, presented just after
(5.4) in \cite{lov2}, shows that 
$$
\pi_1 L f_{i,n} (\alpha, \eta) \;=\; \sum_{z\in \bb Z^d} a(z) \mf
h_{i,n} (\alpha, z) [\eta(z) - \eta(0)]\;.
$$
Since $\eta(z) - \eta(0) = \sum_{1\le j\le d} z_j [\eta(e_j)-\eta(0)]$
for the norm $|\cdot|_\alpha$, the second expression on the right hand
side of (\ref{eq:25}) is equal to
$$
\Big |  \sum_{j=1}^d \{ D_{i,j}(\alpha) - \alpha \sigma_{i,j}
+ h_{i,j,n}( \alpha) \} [\eta(e_j) - \eta(0)]  \Big |_\alpha^2\; ,
$$
where
$$
h_{i,j,n}( \alpha) \;=\; \chi (\alpha) \sum_{z\in\bb Z^d} a(z)\, z_j
\, \mf h_{i,n}(\alpha, \{z\})\;.
$$
By construction, $\mf h_{i,n}(\alpha, \{z\})$ converges to $\mf
f_{i}(\alpha, \{z\})$ uniformly, as $n\uparrow\infty$. In particular,
if we define $D_{i,j}(\alpha)$ as
\begin{equation}
\label{eq:27}
D_{i,j}(\alpha) \;=\; \alpha \sigma_{i,j} \;-\;
\chi (\alpha) \sum_{z\in\bb Z^d} a(z)\, z_j \, \mf f_{i}(\alpha,  \{z\})\; ,
\end{equation}
it not difficult to show from the variational formula for the norm
$|\cdot|_\alpha$ that the second term on the right hand side of
(\ref{eq:25}) also vanishes as $n\uparrow\infty$, uniformly in
$\alpha$. This proves the first statement of the theorem.

Notice that $D_{i,j}(\cdot)$ inherits the smoothness of $\mf
f_{i}(\cdot,  \{z\})$.

It remains to check identity (\ref{eq:26}).  By definition of the
scalar product $\ll \cdot, \cdot \gg_\alpha$ and by identity
(\ref{eq:23}), for any vector $v$ in $\bb R^d$, the left hand side of
(\ref{eq:26}) with the sequence $- \chi(\alpha) f_{i,n} (\alpha, \eta)$
in place of $\mf f_{i,n} (\alpha, \eta)$, is equal to
$$
- \chi(\alpha)^2 < \sum_{j=1}^d v_j \mf h_{j ,n} (\alpha, \cdot), 
\sum_{j=1}^d v_j \mf L_\alpha \mf h_{j ,n} (\alpha, \cdot)>\; .
$$
Since $- \mf L_\alpha \mf h_{j ,n}$ converges to $\mf w_j$ in $\mf
H_{-1}$, uniformly in $\alpha$, and since $\mf h_{j ,n}$ is
$(n,\alpha)$-uniformly bounded in $\mf H_1$, the limit of the previous
expression is equal to the limit of
$$
\chi(\alpha)^2 \sum_{j,k=1}^d v_j v_k <  \mf h_{j ,n} (\alpha, \cdot), 
\mf w_k> \;=\; - \chi(\alpha)^2 \sum_{j,k=1}^d v_j v_k 
\sum_{z\in\bb Z^d_*} \mf h_{j ,n} (\alpha, \{z\}) z_k a(z) \; .
$$
The last identity follows from the explicit formula for $\mf w_k$.
Notice that a factor $1/2$ appeared because in the definition of the
inner product $<\cdot, \cdot>$, there is the term $(n+1)^{-1}$.
By construction, $\mf h_{j ,n} (\alpha, \{z\})$ converges uniformly in
$\alpha$ to $\mf f_{j} (\alpha, \{z\})$. In particular, the previous
sum converges uniformly to
$$
- \chi(\alpha)^2 \sum_{j,k=1}^d v_j v_k 
\sum_{z\in\bb Z^d_*} \mf f_{j} (\alpha, \{z\}) z_k a(z)
$$
By definition (\ref{eq:27}) of the diffusion coefficient
$D_{i,j}(\alpha)$, the previous expression is equal to
$$
\chi(\alpha) \, v \cdot \{ D(\alpha) - \alpha \sigma\} v\; .
$$
This concludes the proof of the theorem. \qed 

We conclude this section with a technical lemma needed above.

\begin{lemma}
\label{s12}
For each $\varepsilon >0$ and $k\ge 1$, there exists $\delta>0$ such
that
\begin{eqnarray*}
\!\!\!\!\!\!\!\!\!\!\!\!\!\!\! &&
\sup_{|\alpha - \beta|\le \delta}
\Vert \mf L_\alpha \mf f_{i,\lambda} (\alpha, \cdot)
- \mf L_\beta \mf f_{i,\lambda} (\beta, \cdot) \Vert_{-1} \;\le\; 
\varepsilon \\
\!\!\!\!\!\!\!\!\!\!\!\!\!\!\! &&
\sup_{|\alpha - \beta|\le \delta} \lambda \sum_{n\ge 1} (1+n)^k
\sum_{A\in \mc E_{*,n}} \{ \mf f_{i,\lambda} (\alpha, A)
- \mf f_{i,\lambda} (\beta, A)\}^2 \;\le\; \varepsilon
\end{eqnarray*}
for all $0<\lambda<1$.
\end{lemma}

The proof of this lemma is implicit in the proof of the regularity of
$f_\lambda (\cdot, A)$ presented in Theorem 5.1 and Lemma 5.2 of
\cite{lov1}. We just need to write the equation for $\mf f_{i,\lambda}
(\alpha, \cdot) - \mf f_{i,\lambda} (\beta, \cdot)$ as a resolvent
equation with a right hand side which is a function in $\mf H_{-1}$
times $O (\alpha - \beta)$. Details are left to the reader.

\section{Technical bounds}
\label{sec4}

We present in this section some technical lemmas and some
computations omitted in section \ref{sec3}. 

\subsection{Computation of $N^2 L_N^* \psi^N_{t,\mf
    f}/\psi^N_{t,\mf f}$.} \label{sec4.0} 
Since $L_N^*$ is the generator of the exclusion process
associated to the transition probability $p^*(y) = p(-y)$, 
$$
\frac{ N^2 L_N^* \psi^N_{t,\mf f} (\eta)} {\psi^N_{t,\mf f}(\eta)} \;=\;
N^2 \sum_{x,y\in\bb T_N^d} \eta(x) [ 1-\eta (x+y)] p^*(y) \Big\{
\frac{\psi^N_{t,\mf f} (\sigma^{x,x+y} \eta)}
{\psi^N_{t,\mf f} (\eta)} - 1 \Big\} \;.
$$
For each fixed bond $(x,y)$, $\psi^N_{t,\mf f} (\sigma^{x,x+y}
\eta)/ \psi^N_{t,\mf f} (\eta)$ is an expression of order $N^{-1}$
because $\mf f_i ( \cdot , \eta)$ is a smooth function for each fixed
configuration $\eta$. We may therefore expand the exponential up to
the second order. The order one term is exactly $N^2 L_N^* \log
\psi^N_{t,\mf f}$ and is responsible for the first two lines of
(\ref{eq:8}) plus a remainder of order $N^{d-1}$. The second order
term is equal to
\begin{eqnarray*}
\!\!\!\!\!\!\!\!\!\!\!\!\! &&
(1/2) \sum_{x,y\in\bb T_N^d} \eta(x) [ 1-\eta (x+y)] p^*(y) 
\bigg\{ N \{\lambda (t,x+y/N) - \lambda(t,x/N)\}  \\
\!\!\!\!\!\!\!\!\!\!\!\!\! && \qquad\qquad\qquad\qquad\qquad\qquad
- \sum_{i=1}^d \sum_{z\in\bb T_N^d} (\partial_{u_i} \lambda)  (t,z/N) 
\nabla_{x,x+y} (A_\ell \mf f_i) (\eta^M (z) ,\tau_{z} \eta )
\bigg\}^2\; .
\end{eqnarray*}
Since $\ell + s_{\mf f} + A \le M$, the gradient $\nabla_{x,x+y}$ acts
either on the first coordinate or on the second but never on both.
$\mf f_i ( \cdot , \eta)$ being a smooth function, the contribution of
the gradient $\nabla_{x,x+y}$ applied on the first coordinate is at
most of order $M^{-d}$.  Since there are $O(M^{d-1})$ boundary sites
$z$ for which $\nabla_{x,x+y} \eta^M (z)$ does not vanish, the total
contribution of the gradient $\nabla_{x,x+y}$ acting on the first
coordinate of $A_\ell \mf f_i$ is of order $M^{-1}$.

We consider now the set of sites $z$ for which the gradient
$\nabla_{x,x+y}$ acts on the second coordinate of $A_\ell \mf f_i$.
In this case, $z$ should be at a distance smaller than $\ell + A$ from
$x$ and we may replace $(\partial_{u_i} \lambda) (t,z/N)$ by
$(\partial_{u_i} \lambda) (t,x/N)$ paying a price of order $\ell^{d+1}
N^{-1}$. At this point, for a fixed $i$, after a change of variables
$z' = z-x$, we may rewrite the sum appearing inside braces in the
previous formula as
$$
(\partial_{u_i} \lambda)  (t,x/N) \tau_x \nabla_{0,y} 
\sum_{z\in\Lambda_{\ell +A} }  \frac 1{|\Lambda_\ell|} 
\sum_{w\in\Lambda_\ell} \mf f_i (\eta^M (z) ,\tau_{z+w} \eta )\; .
$$
Since the summation over $z$ takes place on $\Lambda_{\ell +A}$, we
may replace $\eta^M (z)$ by $\eta^M (0)$ paying a price of order $\ell
/ M$. In this case the previous sum becomes
$$
(\partial_{u_i} \lambda)  (t,x/N) \tau_x \nabla_{0,y} 
\sum_{z\in\bb Z^d} \mf f_i (\eta^M (0) ,\tau_{z} \eta )
\;=\; (\partial_{u_i} \lambda)  (t,x/N) \tau_x \nabla_{0,y}
\Gamma_{\mf f_i (\eta^M (0) , \cdot )}
$$
because the contribution of each fixed $w$ is the same after
replacing $\eta^M (z)$ by $\eta^M (0)$.

To obtain the third line of (\ref{eq:8}) and the correct order of the
remainder, it remains to expand $N \{\lambda (t,x+y/N) -
\lambda(t,x/N)\}$ and to develop the square.

\subsection{Replacement of 
  $L^* (A_\ell \mf f_i) (\eta^M(0), \eta)$ by $L_{\Lambda_\ell}^*
  (A_\ell \mf f_i) (\eta^\ell(0), \eta)$.} \label{sec4.2}
Observe initially that the generator acts either on the first
coordinate or on the second but never on both because we assumed that
$s_{\mf f} + \ell \le M$. Hence, we have to show that the action of
the generator on the first coordinate is negligible. This is the
content of the next result.

\begin{lemma}
\label{s4}
Fix a function $\mf f$ in $\mf F$, a smooth function $G: \bb R_+\times
\bb T^d\to \bb R$ and assume that $M$ satisfies conditions
(\ref{eq:16}). For every $T>0$,
\begin{eqnarray*}
\!\!\!\!\!\!\!\!\!\!\!\!\! &&
\lim_{\ell \to \infty} \limsup_{N\to \infty} \\
\!\!\!\!\!\!\!\!\!\!\!\!\! && \quad
\Big\vert \int_0^T dt \, \int \, N^{1-d} \sum_{z\in\bb T_N^d}
G (t,z/N) \tau_z (L^* - L_{\Lambda_\ell}^*) (A_\ell \mf f)
(\eta^M(0), \eta) \, f_t^N \, d\nu_\alpha^N \Big\vert \;=\; 0\; .
\end{eqnarray*}
\end{lemma}

Notice that in $L_{\Lambda_\ell}^* (A_\ell \mf f) (\eta^M(0), \eta)$,
the generator is acting only on the second coordinate because $\ell\le
M$.

\begin{proof}
  Let $\mf f_1(\alpha, \eta) = (\partial_\alpha \mf f) (\alpha,
  \eta)$. Since $\mf f(\alpha, \cdot)$ is a smooth function, the
  contribution of $(L^* -L_{\Lambda_\ell}^*) (A_\ell \mf f_i)
  (\eta^M(0), \eta)$ is equal to
\begin{equation}
\label{eq:18}
N^{1-d} M^{-d} \sum_{z\in\bb T_N^d} G(t,z/N) \tau_z 
\sum_{\substack{x \in \Lambda_M \\ x+y \not \in \Lambda_M}} 
\eta(x) [ 1-\eta (x+y)] p^*(y) (A_\ell \mf f_1) (\eta^M (0) ,\eta )
\;+\; o_N(1)
\end{equation}
plus a similar term with a negative sign and $x+y$ in $\Lambda_M$, $x$
not in $\Lambda_M$. Here the remainder $o_N(1)$ is of order $N/M^{d+1}$.
From this point, the proof is divided in several steps.

\noindent{\bf Step 1.} 
The first one consists in translating the local functions $\eta(x) [
1-\eta (x+y)]$, which lies at the boundary of $\Lambda_M$, by few
steps in order to have their support contained in $\Lambda_M$. For
this purpose, it is enough to show that for every fixed $y$,
\begin{equation}
\label{eq:17}
N^{1-d} M^{-d} \sum_{z\in\bb T_N^d} G(t,z/N) \tau_z 
\sum_{\substack{x \in \Lambda_M \\ x+y \not \in \Lambda_M}} 
\tau_x W (A_\ell \mf f_1) (\eta^M (0) ,\eta )
\end{equation}
is negligible if $W = h - \tau_{e_1} h$ for some local function $h$.
Here and below, a function $H_{N,\ell}(t, \eta)$ is said to be
negligible if
$$
\lim_{\ell \to \infty}\limsup_{N\to \infty} \Big\vert 
\int_0^T dt \, \int H_{N,\ell} (t,\eta) \, f_t^N \,
d\nu_\alpha^N \Big\vert \;=\; 0
$$
for all $T>0$. Since there exists a finite constant $C_0$ such that
$H_N(\mu^N \,|\, \nu_\alpha^N) \le C_0 N^{d}$ for all measure $\mu^N$,
by the entropy inequality, Feynman-Kac formula and the variational
formula for the largest eigenvalue of a symmetric operator, to prove
that a function is negligible, it is enough to show that 
\begin{equation}
\label{eq:29}
\lim_{\ell \to \infty}\limsup_{N\to \infty}  \int_0^T dt \sup_{f} \Big\{
\int H(t,\eta) \, f \,d\nu_\alpha^N  - \varepsilon N^{2-d} D_N(f)\Big\}
\;\le\; 0
\end{equation}
for every $\varepsilon>0$. Here, the supremum is carried over all
densities $f$ and $D_N(f)$ is the Dirichlet form given by $D_N(f) = <-
L_N \sqrt{f}, \sqrt{f}>$, where $<\cdot, \cdot>$ stands for the inner
product in $L^2(\nu_\alpha^N)$.

Since the local function $W$ has mean zero with respect to all
canonical invariant states, $W=L_\Lambda^s w$ for some finite set
$\Lambda$ and some local function $w$, where $L_\Lambda^s$ stands for
the symmetric part of the generator $L$ restricted to the set
$\Lambda$.  In particular, we need only to show that
$$
N^{1-d} M^{-d} \sum_{z\in\bb T_N^d} G(t,z/N) \tau_z 
\sum_{\substack{x \in \Lambda_M \\ x+y \not \in \Lambda_M}} 
\tau_x (\nabla_b w) (A_\ell \mf f_1) (\eta^M (0) ,\eta )
$$
is negligible for a fixed bond $b=(b_1,b_2)$ and a fixed local
function $w$. Fix $0\le t\le T$, a density $f$ with respect to
$\nu_\alpha^N$ and consider the linear term in variational formula
(\ref{eq:29}):
$$
N^{1-d} M^{-d} \sum_{z\in\bb T_N^d} G(t,z/N) \sum_{\substack{x \in \Lambda_M
    \\ x+y \not \in \Lambda_M}}  \int  \tau_x (\nabla_b w) 
(A_\ell \mf f_1) (\eta^M (0) ,\eta ) \, \tau_{-z} f \,d\nu_\alpha^N\; ,
$$
where we performed a change of variables $\xi = \tau_z \eta$. Since
$\tau_x \nabla_b = \nabla_{b+x} \tau_x$, performing a change of
variables $\xi = \sigma^{b+x} \eta$, we may rewrite the previous
expression as
$$
N^{1-d} M^{-d} \sum_{z\in\bb T_N^d} G(t,z/N) \sum_{\substack{x \in \Lambda_M
    \\ x+y \not \in \Lambda_M}}  \int  \tau_x w
(A_\ell \mf f_1) (\eta^M (0) ,\eta ) \, \nabla_{b+x} \tau_{-z} f 
\,d\nu_\alpha^N
$$
plus a term of order $N M^{-d-1}$. This term appears when taking
the difference $\nabla_{b+x}(A_\ell \mf f_1) (\eta^M (0) ,\eta )$
which is absolutely bounded by $C M^{-d}$. 

Rewrite the difference $ a- b = \tau_{-z} f (\sigma^b \eta) -
\tau_{-z} f (\eta)$ as $(\sqrt{a} - \sqrt{b})(\sqrt{a} + \sqrt{b})$
and apply the elementary inequality $2ab\le \gamma a^2 + \gamma^{-1}
b^2$, which holds for every $\gamma>0$ to estimate the previous
expression by $C \varepsilon^{-1} M^{-2} + \varepsilon N^{2-d}
D_N(f)$. This proves that (\ref{eq:17}) is negligible, concluding the
first step.

\noindent{\bf Step 2.}
Once that all functions have been translated to have its support
contained in $\Lambda_M$, we take advantage of the fact that each
function which appears in (\ref{eq:18}) at one side of the boundary,
appears also at the other side with reversed sign.  In particular,
adding the intermediary terms to complete a telescopic sum, after
(\ref{eq:17}), (\ref{eq:18}) can be rewritten as
$$
N^{1-d} M^{1-d} \sum_{j=1}^m \sum_{z\in\bb T_N^d} G(t,z/N) \tau_z 
\sum_{x \in \Lambda_{M-A}} (\tau_x h_j)(\eta)\,  
(A_\ell \mf f_1) (\eta^M (0) ,\eta )
$$
for a family of local functions $h_j = g_j - \tau_{e_{i}}g_j$ for
some $1\le i\le d$. Here $m$ is a finite integer which depends on
$p(\cdot)$ only. In particular, the local functions $h_j$ have mean
zero with respect to all canonical invariant measures. Here again, $A$
is taken large enough for the support of each local function $\tau_x
h_j$ to be contained in $\Lambda_M$. We claim that such a term is
negligible.

Since all local functions $h$ which have mean zero with respect to all
canonical invariant measures can be expressed as $L^s_{\Lambda} h_0$
for some finite set $\Lambda$ and some local function $h_0$, fix a
bond $b$, a local function $h_0$ and consider the linear term in
(\ref{eq:29}): 
$$
N^{1-d} M^{-d} \sum_{z\in\bb T_N^d} G(t,z/N)  
\sum_{x \in \Lambda_{M-A}} \int \tau_x (\nabla_b h_0)\,  (A_\ell \mf f_1) 
(\eta^M (0) ,\eta ) \, \tau_{-z} f\, d\nu_\alpha^N\; .
$$
Since $\tau_x \nabla_b = \nabla_{b+x} \tau_x$, a change of
variables $\xi = \sigma^{b+x} \eta$, similar to the one performed in
the first part of the proof, permits to write the previous expression
as
\begin{eqnarray}
\label{eq:19}
\!\!\!\!\!\!\!\!\!\!\!\! &&
\frac{N^{1-d}}{M^d} \sum_{z\in\bb T_N^d} G(t,z/N)  
\sum_{x \in \Lambda_{M-A}} \int \, \tau_x h_0 \,  (A_\ell \mf f_1) 
(\eta^M (0) , \sigma^{b+x} \eta ) \, \nabla_{b+x} \tau_{-z} f\,
d\nu_\alpha^N \\
\!\!\!\!\!\!\!\!\!\!\!\! && \quad +\;
\frac{N^{1-d}}{M^d} \sum_{z\in\bb T_N^d} G(t,z/N)  
\sum_{x \in \Lambda_{M-A}} \int \, \tau_x h_0 \,  
\nabla_{b+x} (A_\ell \mf f_1)  (\eta^M (0) , \eta ) \, \tau_{-z} f\,
d\nu_\alpha^N\; . 
\nonumber
\end{eqnarray}
We claim that both terms can be estimated by $ \varepsilon N^{2-d}
D_N(f)$ and an expression which vanishes as $N\uparrow\infty$ and then
$\ell \uparrow\infty$. Notice that in the second term, the gradient
$\nabla_{b+x}$ is acting only on the second coordinate.

Consider the first line of (\ref{eq:19}). Repeating the arguments
presented at the end of the first step, we may bound this integral by
the sum of $\varepsilon N^{2-d} D_N(f)$ and
$$
\frac{C\varepsilon^{-1}}{N^d M^d} \sum_{z\in\bb T_N^d}  
\sum_{x \in \Lambda_{M-A}} \int \,   (A_\ell \mf f_1) 
(\eta^M (0) , \sigma^{b+x} \eta )^2 \, \Big\{ \tau_{-z} f (\eta)
+ \tau_{-z} f (\sigma^{b+x} \eta) \Big\} d\nu_\alpha^N 
$$
for some finite constant $C$. Notice that we got an extra factor
$N^{-1}$ in this passage and that we included $G$ and $h_0$ in the
constant.  We perform a change of variables $\xi = \sigma^{b+x} \eta$
and denote by $\bar f$ the average of the translations of $f$: $\bar f
= N^{-d} \sum_{z\in\bb T_N^d} \tau_z f$ to rewrite the previous sum as
$$
C\varepsilon^{-1} \int \,   (A_\ell \mf f_1) (\eta^M (0) , \eta )^2 \, 
\bar f (\eta) \, d\nu_\alpha^N \; +\; O(\ell^d M^{-d}) \;.
$$
Here we took advantage of the fact that $(A_\ell \mf f_1) (\eta^M
(0) , \sigma^{b+x} \eta ) = (A_\ell \mf f_1) (\eta^M (0) , \eta )$
unless $x$ belongs to $\Lambda _{\ell}$. Since $\mf f_1(\cdot ,
\eta)$ is a smooth function, uniformly in $\eta$, the integral in the
previous expression is less than or equal to
$$
C\varepsilon^{-1}   
\int \,   (A_\ell \mf f_1) (\eta^\ell (0) , \eta )^2 \, 
\bar f (\eta) \, d\nu_\alpha^N \;+\; C\varepsilon^{-1}    
\int \,   \Big\{ \eta^M (0) - \eta^\ell (0) \Big\}^2 \, 
\bar f (\eta) \, d\nu_\alpha^N  \;.
$$
The usual proof of the two blocks estimate permits to show that the
second integral can be estimated by $\varepsilon N^{2-d} D_N(f)$ and
an expression which vanishes as $N\uparrow\infty$ and then $\ell
\uparrow\infty$. We leave the details to the reader. In contrast, the
usual proof of the one block estimate permits to show that the limit,
as $N\uparrow\infty$, of the first integral minus $\varepsilon N^{2-d}
D_N(f)$ is bounded by
$$
C\varepsilon^{-1}  \sup_{K}
\int \,   \Big\{ \frac 1{|\Lambda_{\ell-A}|} \sum_{y\in
  \Lambda_{\ell-A}} \tau_y  \mf f_1 (K/|\Lambda_{\ell}| , \eta ) \Big\} ^2 \, 
d\mu_{\Lambda_{\ell}, K} \; .
$$
In this formula, $\mu_{\Lambda_{\ell}, K}$ stands for the canonical
measure on $\Lambda_\ell$ concentrated on configurations with $K$
particles and the supremum is carried over all integers $0\le K\le
|\Lambda_{\ell}|$. Divide the average in $\Lambda_\ell$ in two
averages and recall from Lemma A.7 in \cite{loy1} that the
Radon-Nikodym derivative $d\mu_{\Lambda_{\ell}, K}/
d\nu^{\Lambda_{\ell}}_{K/|\Lambda_{\ell}|}$ is bounded, uniformly in
$K$, provided $\nu^{\Lambda_{\ell}}_\beta$ stands for the grand
canonical measure on $\Lambda_{\ell}$ with density $\beta$. The
previous expression is thus less than or equal to
$$
C\varepsilon^{-1}  \sup_{0\le\beta\le 1}
\int \,   \Big\{ \frac 1{|\Lambda_{\ell, 1}|} \sum_{y\in
  \Lambda_{\ell, 1}} \tau_y  \mf f_1 (\beta , \eta ) \Big\} ^2 \, 
d\nu^{\Lambda_{\ell}}_{\beta} \; .
$$
In this formula, $\Lambda_{\ell, 1}$ stands for one half of the
cube $\Lambda_\ell$.  Since $\mf f_1 (\alpha , \cdot )$ is local
function, with uniform support and which has mean zero with respect to
$\nu_\alpha^N$, the previous expression is of order $\ell^{-d}$ because
$\nu_\alpha^N$ is a product measure. This conclude the estimation of the
first term in (\ref{eq:19}).

We turn now to the second term of (\ref{eq:19}). Notice that the
gradient $\nabla_{b+x} (A_\ell \mf f_1)$ $(\eta^M (0) , \eta )$
vanishes if $x$ does not belong to $\Lambda_{\ell+A}$. In particular,
\begin{equation}
\label{eq:20}
\sum_{x \in \Lambda_{M-A}}  \tau_x h_0 \,  
\nabla_{b+x} (A_\ell \mf f_1)  (\eta^M (0) , \eta )
\;=\;
\sum_{x \in \Lambda_{\ell+A}} \tau_x h_0 \,  
\nabla_{b+x} (A_\ell \mf f_1)  (\eta^M (0) , \eta )
\end{equation}
is bounded by a constant which does not depend on $N$. On the other
hand, for every $0\le K\le |\Lambda_M|$, repeating the computation
presented in the second paragraph of the second step, from the end to
the beginning, we obtain that
\begin{eqnarray*}
\!\!\!\!\!\!\!\!\!\!\!\! &&
\sum_{x \in \Lambda_{\ell+A}} \int \, \tau_x h_0 \,  
\nabla_{b+x}   (A_\ell \mf f_1)  (\eta^M (0) , \eta )
\, d\mu_{\Lambda_M, K} \\
\!\!\!\!\!\!\!\!\!\!\!\! && \quad \;=\;
\sum_{x \in \Lambda_{\ell+A}} \int \, (\tau_x \nabla_{b} h_0) \,  
(A_\ell \mf f_1)  (\eta^M (0) , \eta )
\, d\mu_{\Lambda_M, K}\; .
\end{eqnarray*}
Summing over all bonds $b$, we recover $L^s h_0 = h = g -
\tau_{e_i}g$, for some local function $g$ and some $1\le i\le d$. The
previous expression is thus equal to
\begin{eqnarray*}
\!\!\!\!\!\!\!\!\!\!\!\! &&
\sum_{x \in \partial_{i}^{-} \Lambda_{\ell+A}} \int \, (\tau_x g) \,  
(A_\ell \mf f_1)  (\eta^M (0) , \eta ) \, d\mu_{\Lambda_M, K} \\
\!\!\!\!\!\!\!\!\!\!\!\! && \quad
\;-\; \sum_{x \in \partial_{i}^{+} \Lambda_{\ell+A}} \int \, (\tau_x g) \,  
(A_\ell \mf f_1)  (\eta^M (0) , \eta ) \, d\mu_{\Lambda_M, K}\; ,
\end{eqnarray*}
where $\partial_{i}^{-} \Lambda_{\ell+A}$ stands for the lower
boundary in the $i$-th direction of $\Lambda_{\ell+A}$ and
$\partial_{i}^{+} \Lambda_{\ell+A}$ for the upper boundary. In
particular, $x$ belongs to $\partial_{i}^{\pm} \Lambda_{\ell+A}$ if it
belongs to $\Lambda_{\ell+A}$ and $\pm x_i =\ell+A$. Since the measure
$\mu_{\Lambda_M, K}$ is uniform,
$$
E_{\mu_{\Lambda_M, K}}[(\tau_x g) g' ] \;=\; E_{\mu_{\Lambda_M, K}}
[(\tau_y g) g']
$$
if the support of $\tau_x g$ and the one of $\tau_y g$ do not
intersect the one of $g'$. Therefore, choosing $A$ large enough, the
previous sum vanishes. This proves that the function (\ref{eq:20}) has
mean zero with respect to all canonical invariant measures.

At this point, we follow the classical approach of nongradient systems
(cf. \cite{kl}, Chapter 7) to estimate the second term of
(\ref{eq:19}) using the standard Rayleigh-Schroedinger perturbation
theorem for the largest eigenvalue of a symmetric operator. After a
few steps we bound the difference of the second term of (\ref{eq:19})
with $\varepsilon N^{2-d} D_N(f)$ by
$$
\frac{N^{2-d} \varepsilon}{M^d} \sum_{z\in\bb T_N^d} \sup_{K}
\Big\{ \frac{G(t,z/N)}{N \varepsilon^2}  \int B \, f \,
d\mu_{\Lambda_M, K} \;-\; <-L^s_{\Lambda_M} \sqrt{f}, 
\sqrt{f}>_{\mu_{\Lambda_M, K}} \Big\}\;.
$$
In this formula, $B$ stands for the function (\ref{eq:20}), the
supremum is carried over all integers $0\le K\le |\Lambda_M|$ and $<
\cdot, \cdot >_{\mu_{\Lambda_M, K}}$ is the inner product in
$L^2(\mu_{\Lambda_M, K})$. Since the spectral gap of the generator of
the symmetric exclusion process in $\Lambda_M$ is of order $M^2$ and
$M^2 N^{-1}$ vanishes as $N\uparrow\infty$, by the perturbation
theorem for the largest eigenvalue of a symmetric operator, the
previous expression is less than or equal to
$$
\frac{C}{M^d \varepsilon^3} \sup_{K}
<(-L^s_{\Lambda_M})^{-1} B , B>_{\mu_{\Lambda_M, K}}\;.
$$
Consider the linear term in the variational formula for the
$H_{-1}$ norm of $B$. It is given by $2 <B , f>_{\mu_{\Lambda_M, K}}$
for some function $f$ in $L^2(\mu_{\Lambda_M, K})$. Since $B$ has mean
zero with respect to all canonical invariant measures, this is in fact
a covariance that we estimate by $C_0(\ell) M^2 + C_1 M^{-2} <f ,
f>_{\mu_{\Lambda_M, K}}$. By the spectral gap for the symmetric
exclusion process, the second term is bounded by $<(-L^s_{\Lambda_M})
f , f>_{\mu_{\Lambda_M, K}}$ if we choose $C_0$ sufficiently small.
Therefore, $<(-L^s_{\Lambda_M})^{-1} B , B>_{\mu_{\Lambda_M, K}}$ is
bounded by $C(\ell) M^2$. Since we are in dimension $d\ge 3$, the
last displayed equation vanishes as $N\uparrow \infty$. This proves
that the second term in (\ref{eq:19}) may be estimated by $\varepsilon
N^{2-d} D_N(f)$ and an expression which vanishes as $N\uparrow\infty$.
\end{proof}

We have just proved that we may replace $L^*$ by $L^*_{\Lambda_\ell}$
in (\ref{eq:8}).  We show now that we can replace the average
$\eta^M(0)$ by the average $\eta^\ell (0)$.

\begin{lemma}
\label{s5}
Fix a function $\mf f$ in $\mf F$, a smooth function $G: \bb R_+\times
\bb T^d\to \bb R$ and assume that $M$ satisfies the
conditions (\ref{eq:16}). For every $T>0$,
\begin{eqnarray*}
\!\!\!\!\!\!\!\!\!\!\!\! &&
\lim_{\ell \to \infty} \limsup_{N\to \infty} \Big\vert \int_0^T dt \, 
\int \, N^{1-d} \sum_{z\in\bb T_N^d} G (t,z/N) \\
\!\!\!\!\!\!\!\!\!\!\!\! && \qquad\qquad\qquad\qquad\quad 
\tau_z \Big\{ L^*_{\Lambda_\ell} (A_\ell \mf f)
(\eta^M(0), \eta) - L^*_{\Lambda_\ell} (A_\ell \mf f)
(\eta^\ell (0), \eta) \Big\} \, f_t^N \, d\nu_\alpha^N \Big\vert \;=\; 0\; .
\end{eqnarray*}
\end{lemma}

\begin{proof}
We have seen in the proof of the previous theorem that it is enough to
show that 
$$
N^{1-d} \sum_{z\in\bb T_N^d}
G (t,z/N) 
\tau_z \Big\{ L^*_{\Lambda_\ell} (A_\ell \mf f)
(\eta^M(0), \eta) - L^*_{\Lambda_\ell} (A_\ell \mf f)
(\eta^\ell (0), \eta) \Big\}
$$
is negligible.

Consider a class of function $B(\beta, \eta)$, $0\le \beta\le 1$,
whose support is contained in $\Lambda_\ell$. Repeating the well known
steps of the proof of the one block estimate we obtain that
$$
\int B(\eta^M(0), \eta) f(\eta) \, d\nu_\alpha^N  \;=\;
\sum_{K=0}^{|\Lambda_M|} C_K(f)
\int B(K/|\Lambda_M|, \eta) f_{M,K} (\eta) \, d\mu_{\Lambda_M,K}\;, 
$$
where, 
$$
C_K(f) \;=\; \int \mb 1\{ \sum_{x\in\Lambda_M} \eta(x) =K\} f \,
d\nu_\alpha^N\;,\quad
f_{M,K} (\eta) \;=\; \frac {f_M}{\int f_{M} (\eta) \, d\mu_{\Lambda_M,K}}
$$
and $f_M$ is the conditional expectation $E_{\nu_\alpha^N}[f \,|\,
\mc F_M]$. Here, for a set $\Lambda$, $\mc F_\Lambda$ stands for the
$\sigma$-algebra generated by $\{\eta (z),\, z\in \Lambda\}$. At this
point, $B(K/|\Lambda_M|, \cdot)$ is a local function with support in
$\Lambda_\ell$ and we repeat the procedure for $f_{M,K}$,
$\mu_{\Lambda_M,K}$ in place of $f$, $\nu_\alpha^N$. We obtain in this
way that the previous sum is equal to
$$
\sum_{K=0}^{|\Lambda_M|} C_K(f) \sum_{k=0}^{|\Lambda_\ell|}
C_k(f_{M,K}) \int B(K/|\Lambda_M|, \eta) f_{M,K,\ell,k} (\eta) 
\, d\mu_{\Lambda_\ell,k}
$$
with the obvious definitions for $C_k(f_{M,K})$, $f_{M,K,\ell,k}$. 

Using that the Dirichlet form is convex, we may estimate
$$
\int \, N^{1-d} \sum_{z\in\bb T_N^d}
G (t,z/N) B(\eta^M(0), \eta)  \, (\tau_{-z} f) \, d\nu_\alpha^N
\;-\; \varepsilon N^{2-d} D_N(f)
$$
by
\begin{eqnarray}
\label{eq:9}
\!\!\!\!\!\!\!\!\!\!\!\! &&
N^{-d} \sum_{z\in\bb T_N^d} \sum_{K=0}^{|\Lambda_M|} C_K(f^z) 
\sum_{k=0}^{|\Lambda_\ell|} C_k(f^z_{M,K}) \\
\!\!\!\!\!\!\!\!\!\!\!\! && \quad
\Big\{
G (t,z/N) N \int B(K/|\Lambda_M|, \eta) f^z_{M,K,\ell,k} (\eta) 
\, d\mu_{\Lambda_\ell,k} - \frac{\varepsilon N^{2}}{|\Lambda_\ell|}
D_{\Lambda_\ell}(f^z_{M,K,\ell,k}, \mu_{\Lambda_\ell,k}) \Big\}\;.
\nonumber
\end{eqnarray}
In this formula, $f^z = \tau_{-z}f$ and $D_{\Lambda_\ell}(\cdot,
\mu_{\Lambda_\ell,k})$ is the Dirichlet form associated to the
generator $L^s_{\Lambda_\ell}$ and the reversible measure
$\mu_{\Lambda_\ell,k}$. Assume that $B(K/|\Lambda_M|, \eta)$ has mean
zero with respect to all invariant states $\mu_{\Lambda_\ell,k}$,
which is the case of the function we are considering in this lemma.
By the Rayleigh-Schroedinger perturbation theorem for the largest
eigenvalue of a symmetric operator, the expression inside braces in
the previous formula is less than or equal to
\begin{equation}
\label{eq:21}
\frac{C |\Lambda_\ell|}{\varepsilon}
< (-L^s_{\Lambda_\ell})^{-1} B(K/|\Lambda_M|, \eta) , 
B(K/|\Lambda_M|, \eta)>_{\mu_{\Lambda_\ell,k}} \;.
\end{equation}

We claim that in the particular case of this lemma, the previous
expression is bounded by $ C \varepsilon^{-1} (K/|\Lambda_M| -
k/|\Lambda_\ell|)^2$.  Indeed, let $h$ be the local function $\mf f
(K/|\Lambda_M|$, $\eta) - \mf f (k/|\Lambda_\ell|, \eta)$. In the case
where $B$ is the function which appears in the statement of the lemma,
the linear term of the variational formula for the $H_{-1}$ norm is
$$
\frac 2{|\Lambda_{\ell'}|} \sum_{y\in \Lambda_{\ell'}} 
\int (L^* \tau_y h)\, f \, d\mu_{\Lambda_\ell,k}\; ,
$$
where $f$ is in $L^2(\mu_{\Lambda_\ell,k})$. Since $L^* \tau_y h$
is a local function which has mean zero with respect to all invariant
measures, we may localize $f$ around $y$, replace the scalar product
by a covariance, use the spectral gap of the symmetric exclusion
process, restricted to a cube whose length depend only on the support
of $h$, and apply Schwarz inequality to bound $<(\nabla_b E[f | \mc
F_\Lambda])^2>$ by $<(\nabla_b f)^2>$. At the end we obtain that the
previous expression is less than or equal to 
$$
\frac C {|\Lambda_{\ell'}|^2} \sum_{y\in \Lambda_{\ell'}} 
<(L^* \tau_y h)^2>_{\mu_{\Lambda_\ell,k}} \;+\;
<- L^s f , f>_{\mu_{\Lambda_\ell,k}} \;.
$$
Since $\mf f(\cdot, \eta)$ is smooth, uniformly in $\eta$, $L^*
\tau_y h$ is absolutely bounded by $|K/|\Lambda_M| -
k/|\Lambda_\ell|\, |$. This proves that (\ref{eq:21}) is bounded above
by $C \varepsilon^{-1} (K/|\Lambda_M| - k/|\Lambda_\ell|)^2$.

Up to this point we proved that the expression inside braces in
(\ref{eq:9}) is bounded above by $C \varepsilon^{-1} (K/|\Lambda_M| -
k/|\Lambda_\ell|)^2$. Recalling the definition of the constants
appearing in (\ref{eq:9}), we have that this sum is in fact
$$
\int \, \frac{C}{\varepsilon N^d} \sum_{z\in\bb T_N^d}
\Big\{ \eta^M(0) - \eta^\ell(0) \Big\}^2   \, (\tau_{-z} f) 
\, d\nu_\alpha^N\; .
$$
It remains to apply the two blocks estimate to conclude the proof.
\end{proof}
\medskip

\subsection{Replacement of $H_{i,j}(\eta)$ by $\sigma_{i,j}
  F(\eta^\ell (0)) + J_{i,j}(\eta^\ell (0))$.} \label{sec4.1} Fix a
smooth function $G : \bb T^d\times \bb R_+\to \bb R$ and two function
$\mf f$, $\mf g$ in $\mf F$. Since the local functions $\mf f(\beta,
\cdot)$ have a common finite support, for each fixed $y$, there exists
a finite integer $A$ such that
$$
\nabla_{0,y} \Gamma_{\mf f (\eta^M(0), \cdot)}
\;=\; \nabla_{0,y} \sum_{z\in \Lambda_A} 
\mf f (\eta^M(0), \tau_z \eta)\; .
$$
Since $\mf f(\cdot , \eta)$ are smooth functions, the difference
between the previous expression and
$$
\nabla_{0,y} \sum_{z\in \Lambda_A} 
\mf f (\eta^\ell (0), \tau_z \eta) 
$$
is absolutely bounded by $C(A, \mf f) \, |\eta^M(0)-\eta^\ell(0)|$,
for some finite constant $C(A, \mf f)$. By the two blocks estimate,
the average over $\bb T_N^d$ of this absolute value is negligible.
After this replacement, the third line of (\ref{eq:8}) is seen to be
composed of three different types of terms:
\begin{eqnarray*}
\!\!\!\!\!\!\!\!\!\!\!\! && 
\sum_{x\in\bb T_N^d} G (t, x/N) \, \tau_x \sum_{y\in \bb Z^d} p^*(y) 
\, y_i \, y_j \, \eta(0) [ 1-\eta (y)] \;, \\
\!\!\!\!\!\!\!\!\!\!\!\! && \quad
\sum_{x\in\bb T_N^d} G (t, x/N) \, \tau_x \sum_{y\in \bb Z^d} p^*(y) 
\,y_i  \,\eta(0) [ 1-\eta (y)] \, \Gamma_{y, \mf f}^{A, \ell} (\eta) \;, \\
\!\!\!\!\!\!\!\!\!\!\!\! && \qquad
\sum_{x\in\bb T_N^d} G (t, x/N) \, \tau_x \sum_{y\in \bb Z^d} p^*(y) 
\, \eta(0) [ 1-\eta (y)] \,  \Gamma_{y, \mf f}^{A, \ell}(\eta) \,  
\Gamma_{y, \mf g}^{A, \ell}(\eta) \;,
\end{eqnarray*}
where, for some function $\mf h$ in $\mf F$,
$$
\Gamma_{y, \mf h}^{A, \ell} (\eta) \;=\; \nabla_{0,y} 
\sum_{z\in \Lambda_A} \mf h (\eta^\ell(0), \tau_z \eta)\;.
$$

By the one block estimate, the first sum can be replaced by 
$$
\sum_{x\in\bb T_N^d} G (t, x/N) \, \sigma_{i,j} \,
F(\eta^\ell(x))\;.
$$

We claim that the second sum is negligible because $\eta (0) [1-\eta
(y)] \Gamma_{y, \mf f}^{A, \ell}$ has mean zero with respect to all
canonical invariant measures. Indeed, repeating the steps of the one
block estimate, we are reduced to estimate
$$
\sup_{K} \int \eta(0) [ 1-\eta (y)] \, \nabla_{0,y}
\sum_{z\in\Lambda_A} \mf f (K/|\Lambda_\ell | , \tau_z \eta) 
\, d\mu_{\Lambda_\ell,K} \;,
$$
where the supremum is carried over all $0\le K\le |\Lambda_\ell|$.
A change of variables $\xi = \sigma^{0,y} \eta$ permits to
rewrite the previous expression as
$$
\sup_{K} \int \big\{ \eta(y) - \eta (0) \big\} \, 
\sum_{z\in\Lambda_A} \mf f (K/|\Lambda_\ell | , \tau_z \eta) 
\, d\mu_{\Lambda_\ell,K} \;.
$$
The integral vanishes for each fixed $K$ because
$\mu_{\Lambda_\ell,K}$ is a uniform measure.

The third type of term requires some notation. For a function $\mf h
(\beta, \eta)$, smooth in the first coordinate and with a common finite
support in the second, let
$$
\tilde {\mf h}(\alpha, \beta) \; =\; E_{\nu_\beta}[\mf h(\alpha,
\eta)]\; .
$$
For $1\le i,j\le d$ and $y$ in $\bb Z^d$, let 
\begin{equation*}
\mf h_y^{i,j} (\beta, \eta) \;=\; \eta(0) [1-\eta(y)]\, \nabla_{0,y} 
\sum_{z\in \Lambda_A} \mf f_i (\beta, \tau_z \eta)\,
\nabla_{0,y} \sum_{z\in \Lambda_A} \mf f_j (\beta, \tau_z \eta)
\;.
\end{equation*}
Notice that $\mf h$ is smooth in the first coordinate and have a
common finite support on the second coordinate. Moreover, an
elementary computation shows that
$$
\sum_{y\in \bb Z^d} p^*(y) \tilde {\mf h}_y^{i,j} (\beta, \beta) \;=\;
2 \sum_{x\in\bb Z^d}  < \mf f_i (\beta, \cdot), (-L^s) \tau_x \mf f_j(\beta,
\cdot)>_\beta \; .
$$
In this formula, $< \cdot, \cdot >_\beta$ stands for the inner product
in $L^2(\nu_\beta)$. Denote the right hand side by $J_{\mf f_i, \mf f_j}
(\beta)$. Lemma \ref{s6} below shows that we may replace in
(\ref{eq:8}) the third type of term by
$$
\sum_{x\in\bb T_N^d} G (t, x/N) \, J_{\mf f_i,\mf f_j} 
(\eta^\ell(x))\; .
$$

Up to this point, we proved that the third line of (\ref{eq:8}) is
equal to 
$$
(1/2) \sum_{i,j=1}^d \sum_{x\in\bb T_N^d} (\partial_{u_i}
\lambda)(t,x/N) (\partial_{u_j} \lambda)(t,x/N) \Big\{
\sigma_{i,j} F(\eta^\ell(x)) + J_{n,i,j} (\eta^\ell(x)) \Big\}
$$
plus a term of order $o(N^d)$, where
$$
J_{n,i,j} (\beta) \;=\;
2 \sum_{x\in\bb Z^d}  < \mf f_{i,n} (\beta, \cdot), 
(-L^s) \tau_x \mf f_{j,n}(\beta, \cdot)>_\beta \; .
$$
Recall the definition of the function $J_{i,j} (\beta)$ given in
(\ref{eq:31}). By Theorem \ref{s7}, with the notation introduced in
section \ref{sec3}, the previous sum is equal to
$$
(1/2) \sum_{i,j=1}^d \sum_{x\in\bb T_N^d} (\partial_{u_i}
\lambda)(t,x/N) (\partial_{u_j} \lambda)(t,x/N) \Big\{
\sigma_{i,j} F(\eta^\ell(x)) + J_{i,j} (\eta^\ell(x)) \Big\}
$$
plus a term of order $o(\mf f, N^d)$. To conclude this subsection, it
remains to prove the next result.

\begin{lemma}
\label{s6}
Fix a function $\mf h (\beta, \eta)$ smooth in the first coordinate
and with finite common support in the second. For positive integers
$\ell$, $m$, let
$$
V_{\ell,m}^{\mf h} (\eta) \;=\;  \Big\vert  \frac 1{|\Lambda_m|}
\sum_{y\in\Lambda_m} \mf h(\eta^\ell(0), \tau_y \eta) 
- \tilde {\mf h}(\eta^\ell(0), \eta^\ell (0))\, \Big\vert\; .
$$
Then,
\begin{equation*}
\lim_{m\to\infty} \limsup_{\ell\to \infty} \sup_{f} \Big\{ 
\int \frac 1{N^d} \sum_{x\in\bb T_N^d} \tau_x  V_{\ell,m}^{\mf h}
(\eta)\, f\,
d\nu_\alpha^N \;-\; \varepsilon N^{2-d} D_N(f) \bigg\} \; =\; 0 
\end{equation*}
for all $\varepsilon>0$. 
\end{lemma}

\begin{proof}
Since $\nu_\gamma$ is the Bernoulli product measure, for each $\beta$,
an elementary computation shows that $(\partial_\gamma) \tilde {\mf h}
(\beta, \gamma) = \sum_{x\in\Lambda} <h(\beta, \eta) ; \eta
(x)>_{\gamma}$, where $<\cdot ; \cdot>_\gamma$ stands for the
covariance with respect to $\nu_\gamma$ and $\Lambda$ for a finite set
which contains the common support of the function $\mf h(\beta,
\cdot)$. In particular, the derivative $(\partial_\gamma \tilde {\mf
  h})(\beta, \gamma)$ is uniformly bounded. Hence,
$$
\Big\vert \tilde {\mf h}(\eta^\ell(0), \eta^\ell (0))
- \tilde {\mf h}(\eta^\ell(0), \eta^m (0))\, \Big\vert
\;\le\; C(\mf h) \big\vert \eta^m(0) - \eta^\ell (0) \big\vert
$$
for some finite constant $C(\mf h)$. It follows from the two blocks
estimate that we may replace $\tilde {\mf h}(\eta^\ell(0), \eta^\ell
(0))$ by $\tilde {\mf h}(\eta^\ell(0), \eta^m (0))$ in the definition
of $V_{\ell,m}^{\mf h}$.

Following the classical proof of the one block estimate, we are
reduced to estimate
$$
\sup_{K} \int \Big\vert  \frac 1{|\Lambda_m|}
\sum_{y\in\Lambda_m} \mf h(K/|\Lambda_\ell |, \tau_y \eta) 
- \tilde {\mf h}(K/|\Lambda_\ell |, \eta^m(0))\, \Big\vert \, 
d\mu_{\Lambda_\ell,K}\;,
$$
where the supremum is carried over all $0\le K\le
|\Lambda_\ell|$. For each fixed $\ell$, denote by $K_\ell$ the
integer which maximizes the previous variational formula. There
exists a subsequence $\ell'$ such that
$K_{\ell'}/|\Lambda_{\ell'}|$ converges to some density $\beta$
in $[0,1]$. In particular, the limsup, as $\ell\uparrow\infty$, of
the previous expression is less than or equal to
$$
\sup_{\beta\in[0,1]} \int \Big\vert  \frac 1{|\Lambda_m|}
\sum_{y\in\Lambda_m} \mf h(\beta, \tau_y \eta) 
- \tilde {\mf h}(\beta , \eta^m(0))\, \Big\vert \, 
d\nu_{\beta}
$$
because the finite marginals of the canonical measure converges
to the grand canonical measures. Since $\tilde  {\mf h}(\beta ,
\cdot)$ is a smooth function,
$$
\tilde {\mf h}(\beta , \eta^m(0)) \; =\; \tilde {\mf h}(\beta ,
\beta) \pm C(\eta^m(0) - \beta ) \; =\;
E_{\nu_{\beta}} [\mf h(\beta , \eta)] \pm C(\eta^m(0) - \beta )\;.
$$
In particular, the previous variational formula is bounded above
by
$$
\sup_{\beta\in[0,1]} \int \Big\vert  \frac 1{|\Lambda_m|}
\sum_{y\in\Lambda_m} \mf h(\beta, \tau_y \eta) 
- E_{\nu_{\beta}} [\mf h(\beta , \eta)] \Big\vert \, d\nu_{\beta}
\;+\; C \sup_{\beta\in[0,1]} \int \big\vert
\eta^m(0) - \beta \, \big\vert \, 
d\nu_{\beta}\;.
$$
This expression vanishes as $m\uparrow\infty$ because $\nu_\beta$
is a product measure and $\mf h(\beta, \cdot)$ are local functions
with a finite common support. This concludes the proof of the lemma.
\end{proof}

\subsection{Estimation of the current.} \label{sec4.3}
Fix $i\le i\le d$ and recall the definition of $V_{i}^{\ell} (\eta)$
given just after (\ref{eq:10}). Let
$$
A_{i, N, \ell, \mf f} (t, \eta) \;=\; N^{1-d} \sum_{x\in\bb T_N^d} 
(\partial_{u_i} \lambda) (t, x/N) \tau_x  
V_{i}^{\ell} (\eta)  \;.
$$
By the nongradient estimates, for every $T\ge 0$,
\begin{eqnarray*}
\!\!\!\!\!\!\!\!\!\!\!\!\! &&
\limsup_{\ell\to\infty} \limsup_{N\to\infty} 
\int_0^T dt\, \int \nu_\alpha^N(d\eta) \, f_t^N (\eta) \,
A_{i, N, \ell , \mf f} (t, \eta) \\
\!\!\!\!\!\!\!\!\!\!\!\!\! && \quad \le\;
C_0 \sup_{\alpha \in [0,1]}
\Big |\! \Big| \! \Big|
w_i^*(\alpha, \eta) + \sum_{1\le j\le d} D_{i,j} (\alpha) [
\eta(e_j) -\eta (0)]  
- L^* \mf f_i(\alpha, \eta) \Big |\! \Big| \! \Big|_\alpha^2
\end{eqnarray*}
for some finite constant $C_0$. Here $\V \cdot \V$ is the norm
introduced at the beginning of section \ref{sec6}.  We refer to
section 6 of \cite{loy1} for the proof. Note that we don't need
in the present context the multiscale analysis of \cite{loy1}.
By Theorem \ref{s7} this expression vanishes if we replace $\mf
f_i$ by $\mf f_{i,n}$ and let $n\uparrow\infty$.

\end{document}